\documentclass[twoside,11pt]{article}
\usepackage{amssymb}
\usepackage{amsthm}
\usepackage{amsmath}
\usepackage[all]{xy}
\usepackage{enumerate}
\usepackage{fancyhdr}
\usepackage{a4wide}

\newtheorem{theorem}{Theorem}[section]
\newtheorem{proposition}[theorem]{Proposition}
\newtheorem{corollary}[theorem]{Corollary}
\newtheorem{lemma}[theorem]{Lemma}
\theoremstyle{definition}
\newtheorem*{notation}{Notation}
\newtheorem*{Beweis}{Proof}
\newtheorem{definition}[theorem]{Definition}
\newtheorem{punto}[theorem]{}
\theoremstyle{remark}
\newtheorem{remark}[theorem]{Remark}
\newtheorem{ex}[theorem]{Example}
\newtheorem{exs}[theorem]{Examples}
\newtheorem{remarks}[theorem]{Remarks}
\CompileMatrices
\swapnumbers
\CompileMatrices
\input{tcilatex}
\begin{document}

\title{Uniformly Flat Semimodules}
\author{\textbf{Jawad Y. Abuhlail}\thanks{%
The author would like to acknowledge the support provided by the Deanship of
Scientific Research (DSR) at King Fahd University of Petroleum $\&$ Minerals
(KFUPM) for funding this work through project No. FT100004.} \\
Department of Mathematics and Statistics\\
King Fahd University of Petroleum $\&$ Minerals\\
abuhlail@kfupm.edu.sa}
\date{}
\maketitle

\begin{abstract}
We revisit the notion of flatness for semimodules over semirings. In
particular, we introduce and study a new notion of \emph{uniformly flat
semimodules} based on the exactness of the tensor functor. We also
investigate the relations between this notion and other notions of flatness
for semimodules in the literature.
\end{abstract}

\section*{Introduction}

\qquad The \emph{homological classification of monoids}, suggested by L. A.
Skornjakov \cite{Sko1969a, Sko1969b}, is still an ongoing project attracting
the attention of many experts in Semigroup Theory and Universal Algebra.
Many papers were devoted to study the category $\mathrm{Act}_{S}$ of right $%
S $\emph{-acts} over a monoid $S$ (a \emph{right }$S$\emph{-act} is a set $A$
with a map $\mu :A\times S\longrightarrow A$ such that $a(st)=(as)t$ and $%
a1_{S}=a$ for all $a\in A$ and $s,t\in S$); for more information see the
encyclopedic manuscript of Kilp et al. \cite{KKM2000}. The philosophy in
several of these papers is to model the theory of modules over rings (\emph{%
e.g.} \cite{AF1974}, \cite{Wis1991}) by studying the interplay between the
(categorical) properties of $\mathrm{Act}_{S}$ and the (algebraic)
properties of $S.$

Another approach to study \emph{Abelian monoids} is to consider them as
\emph{semimodules} over the \emph{semiring} $\mathbb{N}_{0}$ of nonnegative
integers \cite{Go19l99a}. This provides us with a richer structure,
motivates a non-additive version of the theory of modules over rings and
opens the door for developing non-Abelian homological algebra \cite{Ina1997}%
. It is worth mentioning that this approach is supported by the important
role that semirings and semimodules play in emerging areas of research like
\emph{idempotent analysis}, \emph{tropical geometry} and several aspects of
theoretical physics \cite{LM2005}, \cite{KM1997} in addition to many
applications in several branches of mathematics and computer science (\emph{%
e.g.} \cite{GM2008}, \cite{Go19l99a}, \cite{HW1998}).

Although some notions of flatness which are different for $S$-acts (\emph{%
e.g.} \cite[Chapter III]{KKM2000}, \cite{B-F2009} and the papers cited
there) coincide for semimodules as shown by Katsov \cite{Kat2004a}, several
notions of flatness which turn out to be the same for modules are in fact
different for semimodules (\emph{e.g.} \emph{flatness} and \emph{%
mono-flatness }\cite{KN2011}). This results in a rich theory of flatness for
semimodules. In this manuscript, we revisit some of these notions and
introduce a new notion of \emph{uniformly flat semimodules} based on the
exactness of the tensor product functor simulating the classical notion of
flat modules over rings.

The motivation for introducing a new notion of flatness for semimodules can
be understood in light of the following observations: the notions of \emph{%
flat} and $k$\emph{-flat} semimodules introduced in \cite{Alt2004} use
Takahashi's tensor products of semimodules \cite{Tak1982a} which are not the
\emph{natural }tensor products. Among the mains disadvantages of such tensor
products is that the category of semimodules over a commutative semiring is
\emph{not }monoidal and that the tensor functor is \emph{not} left adjoint
to the hom functor as one would expect. In fact, Takahashi's tensor products
solve the universal problem related to such structures in the subcategory of
\emph{cancellative} semimodules, but they fail to provide a universal
solution in the whole category of semimodules (see Section 2 for more
details). Moreover, several results use Takahashi's notion of exact
sequences of semimodules \cite{Tak1981} (see also \cite{Go19l99a}), which we
believe is not natural as well; for more details see the recent manuscript
\cite{Abu}. On the other hand, while the notion of flat semimodules
introduced in \cite{Kat2004a} is quite natural, it does not provide a notion
of \emph{relative flatness} w.r.t. a given family of semimodules which
showed to be important in studying several notions related to \emph{pure
exact sequences} of modules over rings (\emph{e.g.} \cite{Wis1991}). This
motivated us to introduce a new notion of flatness, namely that of \emph{%
uniformly flat semimodules}, using the natural tensor products of
semimodules \cite{Kat1997} and what we believe is a more appropriate notion
of exact sequences of semimodules introduced recently in \cite{Abu}.

This paper is organized as follows. In Section one, we recall some
preliminaries about semirings, semimodules and exact sequences of
semimodules. In Section two, we recall the construction of the natural
tensor products of semimodules over semirings, clarify their connection with
Takahashi's tensor products and study some of their properties. In Section
3, we introduce the notion of \emph{uniformly flat\ semimodules} and
investigate its connection with other notions of flatness in the literature.
We also generalize several results known for modules over rings to
semimodules over semirings.

\section{Preliminaries}

\qquad As pointed out in \cite{KN2011}: \textquotedblleft \emph{when
investigating semirings and their representations, one should undoubtedly
use methods and techniques of both ring and lattice theory as well as
diverse techniques and methods of categorical and universal algebra}%
.\textquotedblright\

For the convenience of the readers who might have different backgrounds, and
to make this manuscript as much self-contained as possible, we collect in
this section some definitions, remarks and results that will be used in the
sequel. For unexplained terminology, our main references are \cite{Mac1998}
for Category Theory, \cite{Gra2008} for Universal Algebra and \cite{Wis1991}
for Ring and Module Theory.

\subsection*{Semirings and Semimodules}

\qquad \emph{Semirings} (\emph{semimodules}) are roughly speaking, rings
(modules) without subtraction. Recall that a semigroup $(S,\ast )$ is said
to be \emph{cancellative} iff for any $s_{1},s_{2},s\in S$ we have%
\begin{equation*}
\lbrack s_{1}\ast s=s_{2}\ast s\Longrightarrow s_{1}=s_{2}]\text{ and }%
[s\ast s_{1}=s\ast s_{2}\Longrightarrow s_{1}=s_{2}].
\end{equation*}

\begin{definition}
A \emph{semiring} is an algebraic structure $(S,+,\cdot ,0,1)$ consisting of
a non-empty set $S$ with two binary operations \textquotedblleft $+$%
\textquotedblright\ (addition) and \textquotedblleft $\cdot $%
\textquotedblright\ (multiplication) satisfying the following conditions:

\begin{enumerate}
\item $(S,+,0)$ is an Abelian monoid with neutral element $0_{S};$

\item $(S,\cdot ,1)$ is a monoid with neutral element $1;$

\item $x\cdot (y+z)=x\cdot y+x\cdot z$ and $(y+z)\cdot x=y\cdot x+z\cdot x$
for all $x,y,z\in S;$

\item $0\cdot s=0=s\cdot 0$ for every $s\in S$ (\emph{i.e.} $0$ is \emph{%
absorbing}).
\end{enumerate}
\end{definition}

\begin{punto}
Let $S,S^{\prime }$ be semirings. A map $f:S\rightarrow S^{\prime }$ is said
to be a \emph{morphism of semirings} iff for all $s_{1},s_{2}\in S:$%
\begin{equation*}
f(s_{1}+s_{2})=f(s_{1})+f(s_{2}),\text{ }f(s_{1}s_{2})=f(s_{1})f(s_{2}),%
\text{ }f(0_{S})=0_{S^{\prime }}\text{ and }f(1_{S})=1_{S^{\prime }}.
\end{equation*}
\end{punto}

\begin{punto}
Let $(S,+,\cdot )$ be a semiring. We say that $S$ is

\emph{cancellative} iff the additive semigroup $(S,+)$ is cancellative;

\emph{commutative} iff the multiplicative semigroup $(S,\cdot )$ is
commutative;

\emph{semifield} iff $(S\backslash \{0\},\cdot ,1)$ is a commutative group.
\end{punto}

\begin{exs}
Rings are indeed semirings. The first natural example of a \emph{commutative
}semiring is $(\mathbb{N}_{0},+,\cdot ),$ the set of nonnegative integers.
The semirings $(\mathbb{R}_{0}^{+},+,\cdot )$ and $(\mathbb{Q}%
_{0}^{+},+,\cdot )$ are indeed semifields. Moreover, for any ring $R$ we
have a semiring structure $(,+,\cdot )$ on the set $\mathrm{Ideal}(R)$ of
ideals of $R$ with the usual addition and multiplication of ideals of $R.$
For more examples, the reader may refer to \cite{Go19l99a}.
\end{exs}

\begin{definition}
Let $S$ be a semiring. A \emph{right }$S$\emph{-semimodule} is an algebraic
structure $(M,+,0_{M})$ consisting of a non-empty set $M,$ a binary
operation \textquotedblleft $+$\textquotedblright\ along with a right $S$%
-action%
\begin{equation*}
M\times S\longrightarrow M,\text{ }(m,s)\mapsto ms,
\end{equation*}%
such that:

\begin{enumerate}
\item $(M,+,0_{M})$ is an Abelian monoid with neutral element $0_{M};$

\item $(ms)s^{\prime }=m(ss^{\prime }),$ $(m+m^{\prime })s=ms+m^{\prime }s$
and $m(s+s^{\prime })=ms+ms^{\prime }$ for all $s,s^{\prime }\in S$ and $%
m,m^{\prime }\in M;$

\item $m1_{S}=m$ and $m0_{S}=0_{M}=0_{M}s$ for all $m\in M$ and $s\in S.$
\end{enumerate}
\end{definition}

\begin{punto}
\begin{enumerate}
\item Let $M$ and $M^{\prime }$ be right $S$-semimodules. A map $%
f:M\rightarrow M^{\prime }$ is said to be a \emph{morphism of }$S$\emph{%
-semimodules} (or $S$\emph{-linear}) iff for all $m_{1},m_{2}\in M$ and $%
s\in S:$%
\begin{equation*}
f(m_{1}+m_{2})=f(m_{1})+f(m_{2})\text{ and }f(ms)=f(m)s.
\end{equation*}%
The set $\mathrm{Hom}_{S}(M,M^{\prime })$ of $S$-linear maps from $M$ to $%
M^{\prime }$ is clearly an Abelian monoid under addition. The category of
right $S$-semimodules is denoted by $\mathbb{S}_{S}.$ Analogously, one can
define the category $_{S}\mathbb{S}$ of left $S$-semimodules. A right $S$%
-semimodule $M_{S}$ is said to be \emph{cancellative} iff the semigroup $%
(M,+)$ is cancellative. With $\mathbb{CS}_{S}\subseteq \mathbb{S}_{S}$
(resp. $_{S}\mathbb{CS}\subseteq $ $_{S}\mathbb{S})$ we denote the full
subcategory of cancellative right (left) $S$-semimodules.
\end{enumerate}
\end{punto}

\begin{punto}
Let $S$ be a semiring and $M$ a right $S$-semimodule. A non-empty subset $%
L\subseteq M$ is said to be an $S$\emph{-subsemimodule, }and we write $L\leq
_{S}M,$ iff $L$ is closed under \textquotedblleft $+_{M}$\textquotedblright\
and $ls\in L$ for all $l\in L$ and $s\in S.$
\end{punto}

\begin{ex}
Every Abelian monoid $(M,+,0_{M})$ is an $\mathbb{N}_{0}$-semimodule in the
obvious way. Moreover, the categories $\mathbf{AbMon}$ of Abelian monoids
and the category $\mathbb{S}_{\mathbb{N}_{0}}$ of $\mathbb{N}_{0}$%
-semimodules are isomorphic.
\end{ex}

\begin{punto}
Let $S,T$ be semirings, $M$ a left $S$-semimodule and a right $T$%
-semimodule. We say that $M$ is an $(S,T)$\emph{-bisemimodule} iff $%
(sm)t=s(mt)$ for all $s\in S,$ $m\in M$ and $t\in T.$ For $(S,T)$%
-bisemimodules $M,$ $M^{\prime },$ we call an $S$-linear $T$-linear map $%
f:M\rightarrow N^{\prime }$ a \emph{morphism of }$(S,T)$\emph{- bisemimodules%
} (or $(S,T)$\emph{-bilinear}). The set $\mathrm{Hom}_{(S,T)}(M,M^{\prime })$
of $(S,T)$-bilinear maps from $M$ to $M^{\prime }$ is clearly an Abelian
monoid under addition. The category of $(S,T)$-bisemimodules will be denoted
by $_{S}\mathbb{S}_{T}.$
\end{punto}

Throughout, and unless otherwise explicitly specified, $S$ is an associative
semiring with $1_{S}\neq 0_{S}.$ We mean by an $S$-semimodule a \emph{right}
$S$-semimodule unless something different is mentioned explicitly.

\begin{punto}
\label{zuka}Let $M$ be an $S$-semimodule. An equivalent relation $\equiv $
on $M$ is said to be an $S$\emph{-congruence} iff for any $m,m^{\prime
},m_{1},m_{1}^{\prime },m_{2},m_{2}^{\prime }\in M$ and $s\in S$ we have%
\begin{equation*}
\lbrack (m_{1}\equiv m_{1}^{\prime }\text{ }\wedge \text{ }m_{2}\equiv
m_{2}^{\prime })\Rightarrow m_{1}+m_{2}\equiv m_{1}^{\prime }+m_{2}^{\prime
}]\text{ and }[m\equiv m^{\prime }\Rightarrow ms\equiv m^{\prime }s].
\end{equation*}%
Every $S$-subsemimodule $L\leq _{S}M$ induces two $S$-congruences on $M$
given by%
\begin{equation*}
m_{1} \equiv_{L} \, m_{2}\Leftrightarrow m_{1}+l_{1}=m_{2}+l_{2}\text{ for
some }l_{1},l_{2}\in L;
\end{equation*}
\begin{equation*}
m_{2} \equiv_{[L]} \, m_{2}\Longleftrightarrow m_{1}+l_{1}+m^{\prime \prime
}=m_{2}+l_{2}+m^{\prime \prime }\text{ for some }l_{1},l_{2}\in L\text{ and }%
m^{\prime \prime }\in M.
\end{equation*}
We call the $S$-semimodule $M/L:=M/_{\equiv _{L}}$ the \emph{quotient }(%
\emph{factor})\emph{\ semimodule of }$M$\emph{\ by }$L.$ If $M$ is
cancellative, then $L$ and $M/L$ are cancellative. On the other hand, $%
M/_{\equiv _{\lbrack L]}}$ is obviously cancellative.
\end{punto}

\begin{punto}
Let $M$ be an $S$-semimodule and recall the $S$-congruence relation $\equiv
_{\lbrack 0]}$ on $M$ defined by%
\begin{equation*}
m\equiv _{\lbrack 0]}\text{ }m^{\prime }\text{ }\Longleftrightarrow
m+m^{\prime \prime }=m^{\prime }+m^{\prime \prime }\text{ for some }%
m^{\prime \prime }\in M.
\end{equation*}%
The quotient $S$-semimodule $M/\sim $ is indeed cancellative and we have a
canonical surjection $\mathfrak{c}_{M}:M\longrightarrow \mathfrak{c}(M)$ with%
\begin{equation*}
\mathrm{Ker}(\mathfrak{c}_{M})=\{m\in M\mid m+m^{\prime \prime }=m^{\prime
\prime }\text{ for some }m^{\prime \prime }\in M\}.
\end{equation*}%
The class of cancellative right $S$-semimodules is a \emph{reflective}
subcategory of $\mathbb{S}_{S}$ in the sense that the functor%
\begin{equation*}
\mathfrak{c}:\mathbb{S}_{S}\longrightarrow \mathbb{CS}_{S},\text{ }M\mapsto
M/\equiv _{\lbrack 0]}
\end{equation*}%
is left adjoint to the embedding functor $\mathbb{CS}_{S}\hookrightarrow
\mathbb{S}_{S},$ \emph{i.e.} for any $S$-semimodule $M$ and any \emph{%
cancellative} $S$-semimodule $N$ we have a natural isomorphism of Abelian
monoids $\mathrm{Hom}_{S}(\mathfrak{c}(M),N)\simeq \mathrm{Hom}_{S}(M,N)$
\cite[Page 517]{Tak1981}.
\end{punto}

\begin{proposition}
The category $\mathbb{S}_{S}$ and its full subcategory $\mathbb{CS}_{S}$
have kernels and cokernels, where for any morphism of $S$-semimodules $%
f:M\rightarrow N$ we have%
\begin{equation*}
\mathrm{Ker}(f)=\{m\in M\mid f(m)=0\}\text{ and }\mathrm{\mathrm{Coker}}%
(f)=N/f(M).
\end{equation*}
\end{proposition}

\begin{punto}
We call a subset $Y\subseteq N$ \emph{subtractive} iff $Y=\overline{Y},$ the
\emph{subtractive closure} of $Y,$ where%
\begin{equation*}
\overline{Y}=\{n\in N\mid n+y_{1}=y_{2}\text{ for some }y_{1},y_{2}\in Y\}.
\end{equation*}%
An $S$-semimodule $M$ is said to be \emph{completely subtractive} iff every $%
S$-subsemimodule of $M$ is subtractive.

We call a morphism of $S$-semimodules $\gamma :M\longrightarrow N:$

$k$\emph{-uniform} iff for any $m_{1},m_{2}\in M:$%
\begin{equation}
\gamma (m_{1})=\gamma (m_{2})\Longrightarrow \text{ }\exists \text{ }%
k_{1},k_{2}\in \mathrm{Ker}(\gamma )\text{ s.t. }m_{1}+k_{1}=m_{2}+k_{2};
\label{k-steady}
\end{equation}

$i$\emph{-uniform} iff $\gamma (M)=\overline{\gamma (M)};$

\emph{uniform }iff $\gamma $ is $k$-uniform and $i$-uniform.
\end{punto}

\begin{remark}
The uniform ($k$-uniform, $i$-uniform) morphisms of semimodules were called
\emph{regular} ($k$\emph{-regular, }$i$-\emph{regular}) by Takahashi \cite%
{Tak1982c}. We think that our terminology avoids confusion sine a regular
monomorphism (regular epimorphism) has a different well-established meaning
in the language of Category Theory.
\end{remark}

\begin{punto}
Let $M$ be an $S$-semimodule, $L\leq _{S}M$ an $S$-subsemimodule and
consider the factor semimodule $M/L.$ Then we have a surjective \emph{uniform%
} morphism of $S$-semimodules%
\begin{equation*}
\pi _{L}:=M\rightarrow M/L,\text{ }m\mapsto \lbrack m]
\end{equation*}%
with%
\begin{equation*}
\mathrm{Ker}(\pi _{L})=\{m\in M\mid m+l_{1}=l_{2}\text{ for some }%
l_{1},l_{2}\in L\}=\overline{L};
\end{equation*}%
in particular, $L=\mathrm{Ker}(\pi _{L})$ if and only if $L\subseteq M$ is
subtractive.
\end{punto}

In \cite{Abu} we introduced a new notion of \emph{exact sequences} of
semimodules. Takahashi's exact sequences \cite{Tak1981} shall be called
\emph{semi-exact} in the sequel:

\begin{definition}
We call a sequence of $S$-semimodules%
\begin{equation}
L\overset{f}{\longrightarrow }M\overset{g}{\longrightarrow }N  \label{ABC}
\end{equation}%
\emph{exact }iff $f(L)=\mathrm{Ker}(g)$ and $g$ is $k$-uniform. An exact
sequence%
\begin{equation}
0\longrightarrow L\overset{f}{\longrightarrow }M\overset{g}{\longrightarrow }%
N\longrightarrow 0  \label{ses}
\end{equation}%
is called a \emph{short exact sequence}.
\end{definition}

\begin{punto}
\label{def-exact}(\cite{Abu}) We call a sequence of $S$-semimodules $L%
\overset{f}{\rightarrow }M\overset{g}{\rightarrow }N:$

\emph{proper-exact} iff $f(L)=\mathrm{Ker}(g);$

\emph{semi-exact} iff $\overline{f(L)}=\mathrm{Ker}(g);$

\emph{quasi-exact} iff $\overline{f(L)}=\mathrm{Ker}(g)$ and $g$ is $k$%
-uniform.
\end{punto}

\begin{punto}
We call a (possibly infinite) sequence of $S$-semimodules
\begin{equation}
\cdots \rightarrow M_{i-1}\overset{f_{i-1}}{\rightarrow }M_{i}\overset{f_{i}}%
{\rightarrow }M_{i+1}\overset{f_{i+1}}{\rightarrow }M_{i+2}\rightarrow \cdots
\label{chain}
\end{equation}

\emph{chain complex} iff $f_{j+1}\circ f_{j}=0$ for every $j;$

\emph{exact} (resp. \emph{proper-exact}, \emph{semi-exact, quasi-exact}) iff
each partial sequence with three terms $M_{j}\overset{f_{j}}{\rightarrow }%
M_{j+1}\overset{f_{j+1}}{\rightarrow }M_{j+2}$ is exact (resp. proper-exact,
semi-exact, quasi-exact);
\end{punto}

\begin{punto}
An $S$-semimodule $N$ is said to be a \emph{retract} of an $S$-semimodule $M$
iff there exist a (surjective) $S$-linear map $\theta :M\longrightarrow N$
and an (injective) $S$-linear map $\psi :N\longrightarrow M$ such that $%
\theta \circ \psi =\mathrm{id}_{N}$ (equivalently, $N\simeq \alpha (M)$ for
some \emph{idempotent} endomorphism $\alpha \in \mathrm{End}(M_{S})$). On
the other hand, $N$ is a direct summand of $M$ (\emph{i.e.} $M=N\oplus
N^{\prime }$ for some $S$-subsemimodule $N^{\prime }$ of $M$) if and only if
there exists $\alpha \in \mathrm{Comp}(\mathrm{End}(M_{S}))$ s.t. $\alpha
(M)=N$ where for any semiring $T$ we set%
\begin{equation*}
\mathrm{Comp}(T)=\{t\in T\mid \text{ }\exists \text{ }\widetilde{t}\in T%
\text{ with }t+\widetilde{t}=\mathrm{id}_{T}\text{ and }t\widetilde{t}=0_{T}=%
\widetilde{t}t\}.
\end{equation*}%
Indeed, every direct summand of $M$ is a retract of $M;$ the converse is not
true in general (cf. \cite[Proposition 16.6]{Go19l99a}).
\end{punto}

\begin{lemma}
\label{exact}\emph{(\cite[Proposition 3.10, Corollary 3.11]{Abu}) }Let $A,B$
and $C$ be $S$-semimodules.

\begin{enumerate}
\item $0\longrightarrow A\overset{f}{\longrightarrow }B$ is exact if and
only if $f$ is injective.

\item $B\overset{g}{\longrightarrow }C\longrightarrow 0$ is exact if and
only if $g$ is surjective.

\item $0\longrightarrow A\overset{f}{\longrightarrow }B\overset{g}{%
\longrightarrow }C$ is semi-exact and $f$ is uniform if and only if $A=%
\mathrm{Ker}(g).$

\item $A\overset{f}{\longrightarrow }B\overset{g}{\longrightarrow }%
C\longrightarrow 0$ is semi-exact and $g$ is uniform if and only if $C=%
\mathrm{Coker}(f).$

\item $0\longrightarrow A\overset{f}{\longrightarrow }B\overset{g}{%
\longrightarrow }C\longrightarrow 0$ is exact if and only if $A=\mathrm{Ker}%
(g)$ and $C=\mathrm{Coker}(f).$
\end{enumerate}
\end{lemma}

\qquad The following technical result follows immediately from the
definitions and \cite[Lemmas 1.11, 1.15]{Tak1983}.

\begin{lemma}
\label{comm}

\begin{enumerate}
\item Consider a commutative diagram of $S$-semimodules with $\pi \circ
\iota =\mathrm{id}_{N}$ and $\pi ^{\prime }\circ \iota ^{\prime }=\mathrm{id}%
_{N^{\prime }}:$%
\begin{equation*}
\xymatrix{N \ar[r]^{\tilde{\gamma}} \ar[d]_{\iota} & N' \ar[d]^{\iota '}\\ M
\ar[r]^{\gamma} \ar[d]_{\pi} & M' \ar[d]^{{\pi}'} \\ N
\ar[r]^{\tilde{\gamma}} & N'}
\end{equation*}%
If $\gamma $ is uniform \emph{(}resp. $k$-uniform, $i$-uniform\emph{)}, then
$\widetilde{\gamma }$ is uniform (resp. $k$-uniform, $i$-uniform).

\item Consider a commutative diagram of $S$-semimodules with $\pi \circ i=%
\mathrm{id}_{N},$ $\pi ^{\prime }\circ i^{\prime }=\mathrm{id}_{N^{\prime }}$
and $\pi ^{\prime \prime }\circ i^{\prime \prime }=\mathrm{id}_{N^{\prime
\prime }}:$%
\begin{equation*}
\xymatrix{N \ar[r]^{\tilde{f}} \ar[d]_{\iota} & N' \ar[d]^{\iota '}
\ar[r]^{\tilde{g}} & N'' \ar[d]^{\iota ''}\\ M \ar[r]^{f} \ar[d]_{\pi} & M'
\ar[d]^{{\pi}'} \ar[r]^{g} & M'' \ar[d]^{{\pi}''}\\ N \ar[r]^{\tilde{f}} &
N' \ar[r]^{\tilde{g}} & N''}
\end{equation*}%
If $M\overset{f}{\longrightarrow }M^{\prime }\overset{g}{\longrightarrow }%
M^{\prime \prime }$ is exact \emph{(}resp. proper-exact, semi-exact,
quasi-exact\emph{)}, then $N\overset{\widetilde{f}}{\longrightarrow }%
N^{\prime }\overset{\widetilde{g}}{\longrightarrow }N^{\prime \prime }$ is
exact \emph{(}resp. proper-exact, semi-exact, quasi-exact\emph{)}.
\end{enumerate}
\end{lemma}

Some \emph{redundant} assumptions in \cite[Lemma 4.5]{Abu} do not hold in
some situations which we will handle in this paper. A slight adjustment of
the proof of the above mentioned result yields

\begin{lemma}
\label{diagram}Consider the following commutative diagram of $S$-semimodules%
\begin{equation*}
\xymatrix{L_1 \ar[r]^{f_1} \ar[d]_{\alpha_1} & M_1 \ar[r]^{g_1}
\ar[d]_{\alpha_2} & N_1 \ar[d]_{\alpha_3} \\ L_2 \ar[r]^{f_2} & M_2
\ar[r]^{g_2} & N_2}
\end{equation*}

\begin{enumerate}
\item Let the second sequence be quasi-exact \emph{(i.e.} $\overline{%
f_{2}(L_{2})}=\mathrm{Ker}(g_{2})$ and $g_{2}$ is $k$-uniform\emph{)} and $%
g_{1},$ $\alpha _{1}$ be surjective.

\begin{enumerate}
\item Let $g_{1}\circ f_{1}=0.$ If $\alpha _{2}$ is injective, then $\alpha
_{3}$ is injective.

\item If $\alpha _{3}$ is surjective \emph{(}and $\alpha _{2}$ is $i$-uniform%
\emph{), }then $\alpha _{2}$ is a semi-epimorphism \emph{(}surjective\emph{)}%
.
\end{enumerate}

\item Let the first row be semi-exact \emph{(i.e.} $\overline{f_{1}(L_{1})}=%
\mathrm{Ker}(g_{1})$\emph{) }and $f_{2}$ be injective.

\begin{enumerate}
\item Let $f_{1},$ $\alpha _{2}$ be cancellative and $g_{1}$ be $k$-uniform.
If $\alpha _{1},$ $\alpha _{3}$ are injective, then $\alpha _{2}$ is
injective.

\item Let $g_{2}\circ f_{2}=0.$ If $\mathrm{Ker}(\alpha _{3})=0,$ and $%
\alpha _{2}$ is surjective, then $\alpha _{1}$ is a semi-epimorphism. If
moreover, $\alpha _{1}$ \emph{or} $f_{1}$ is $i$-uniform, then $\alpha _{1}$
is surjective.
\end{enumerate}
\end{enumerate}
\end{lemma}

\section{Tensor products of semimodules}

\qquad Tensor products of semimodules were introduced and investigated by
Takahashi \cite{Tak1981}. However, they did not provide a solution to the
universal problem related to such structures in the whole category of
semimodules. On the other hand, Katsov \cite{Kat1997} considered a different
tensor product in the category of semimodules (over a commutative semiring)
which solved several of the problems that Takahashi's tensor products had.
It is worth mentioning, as Katsov pointed out, that his construction of the
tensor product and the elementary results related to it seem to be folklore (%
\emph{e.g.} Grillet \cite{Gril1969} gave an explicit construction of a \emph{%
non-associative} tensor product in the variety $\mathbf{Sgr}$ of semigroups
and suggested that the same construction works for all algebraic varieties
of Universal Algebra). Varieties in which the tensor products behave nicely
were considered by F. Linton \cite{Lin1966} (see also \cite[Theorem 3.10.3]%
{Bor1994b}).

\subsection*{Construction of tensor products}

\qquad As before, $S$ denotes an associative semiring with $1_{S}\neq 0_{S}.$
With $_{S}\mathbb{S}$ and $\mathbb{S}_{S}$ we denoted the categories of left
and right $S$-semimodule, respectively. For the convention of the reader, we
recall the construction of tensor products of semimodules and some of its
properties (\emph{e.g.} \cite{Kat1997}, \cite{Kat2004a}, \cite{KN2011}):

\begin{punto}
Let $M_{S}$ be a right $S$-semimodule and $_{S}N$ a left $S$-semimodule. An $%
S$-balanced map $g:M\times N\rightarrow G,$ where $G$ is an Abelian monoid,
is a bilinear map such that $g(ms,n)=g(m,sn)$ for all $m\in M,$ $s\in S$ and
$n\in N.$ Let $F$ be the free Abelian monoid with basis $M\times N.$ Every
element of $F$ can be written uniquely as a linear combination of elements
of the set $\{\delta _{(m,n)}\mid (m,n)\in M\times N\}$ where $\delta
_{(m,n)}$ is the Kronecker delta function. Let $\sigma \subseteq F\times F$
be the congruence relation \emph{generated} by the set of all ordered pairs
\begin{equation*}
\{(\delta _{(m_{1}+m_{2},n)},\delta _{(m_{1},n)}+\delta
_{(m_{2},n)}),(\delta _{(m,n_{1}+n_{2})},\delta _{(m,n_{1})}+\delta
_{(m,n_{2})}),(\delta _{(ms,n)},\delta _{(m,sn)})\},
\end{equation*}%
where $m_{1},m_{2},m\in M,$ $n_{1},n_{2},n\in N,$ $s\in S$ and consider
canonical maps%
\begin{equation*}
\pi _{\sigma }:F\longrightarrow F/\sigma \text{ and }\tau :=\pi _{\sigma
}\circ \iota :M\times N\rightarrow F/\sigma .
\end{equation*}%
Let $G$ be an Abelian monoid and $\beta :M\times N\rightarrow G$ an $S$%
-balanced map. Since $F$ is free over $M\times N,$ the map $\beta $ induces
a unique map $\beta ^{\prime }:F\rightarrow G.$ Since $\mathrm{Ker}(\pi
_{\sigma })=\sigma \subseteq \mathrm{Ker}(g),$ there exists a unique map $%
\gamma :F/\sigma \rightarrow G$ such that $\gamma \circ \pi _{\sigma }=\beta
^{\prime }$ (given by $\gamma (\overline{f})=\beta ^{\prime }(f),$ for every
$f\in F$) and so $\gamma \circ \tau =\gamma \circ \pi _{\sigma }\circ \iota
=\beta ^{\prime }\circ \iota =\beta :$%
\begin{equation}
\begin{array}{ccc}
\xymatrix{& M \times N \ar@{^(->}^{\iota}[d] \ar@{.>}[ddl]_{\tau}
\ar^{\beta}[ddr] & \\ & F \ar_{\beta '}[dr] \ar^{\pi_{\sigma}}[dl] & \\
F/\sigma \ar@{.>}_{\gamma}[rr] & & G} &  & \xymatrix{ &M \times N
\ar_{\tau}[ddl] \ar^{\beta}[ddr] & \\ & & \\ M \otimes_{S} N
\ar@{.>}_{\gamma}[rr] & & G}%
\end{array}
\label{uni-ten}
\end{equation}%
So, $M\otimes _{S}N:=(F/\sigma ,\tau )$ is solution for the following
universal problem: For every Abelian monoid $G$ with an $S$-balanced map $%
\beta :M\times N\rightarrow G,$ there exists a unique morphism of monoids $%
\gamma :M\otimes _{S}N\rightarrow G$ that completes the right triangle in (%
\ref{uni-ten}) commutatively.
\end{punto}

\begin{punto}
Let $M_{S}$ a right $S$-semimodule, $_{S}N$ a left $S$-semimodule and $F$
the free Abelian monoid with basis $M\times N.$ Let $N(M)\subseteq F\times F$
be the \emph{symmetric} $S$-subsemimodule generated by the set of elements
of the form%
\begin{equation*}
\begin{array}{ccc}
(\delta _{(m_{1}+m_{2},n)},\delta _{(m_{1},n)}+\delta _{(m_{2},n)}), &  &
(\delta _{(m_{1},n)}+\delta _{(m_{2},n)},\delta _{(m_{1}+m_{2},n)}), \\
(\delta _{(m,n_{1}+n_{2})},\delta _{(m,n_{1})}+\delta _{(m,n_{2})}), &  &
(\delta _{(m,n_{1})}+\delta _{(m,n_{2})},\delta _{(m,n_{1}+n_{2})}), \\
(\delta _{(ms,n)},\delta _{(m,sn)}), &  & (\delta _{(m,sn)},\delta
_{(ms,n)}),%
\end{array}%
\end{equation*}%
and consider the $S$-congruence relation on $F$ defined by%
\begin{equation*}
f\rho g\Longleftrightarrow f+h=g+h^{\prime }\text{ for some }(h,h^{\prime
})\in N(M).
\end{equation*}%
\emph{Takahashi's tensor product} of $M$ and $N$ is defined as $M\boxtimes
_{S}N=F/\rho .$ Notice that there is an $S$-balanced map%
\begin{equation*}
\widetilde{\tau }:M\times N\longrightarrow M\boxtimes _{S}N,\text{ }%
(m,n)\mapsto m\boxtimes _{S}n:=(m,n)/\rho
\end{equation*}%
with the following universal property \cite{Tak1982a}: for every \emph{%
cancellative }Abelian monoid $\widetilde{G}$ and every $S$-balanced map $%
\widetilde{\beta }:M\times N\longrightarrow \widetilde{G}$ there exists a
\emph{unique} morphism of monoids $\widetilde{\gamma }:M\boxtimes
_{S}N\longrightarrow \widetilde{G}$ such that $\widetilde{\gamma }\circ
\widetilde{\tau }=\widetilde{G}.$

The above mentioned property means that $-\boxtimes _{S}-$ plays the role of
a \emph{tensor product} w.r.t. cancellative semimodules. On the other hand,
notice that for every Abelian monoid $G,$ we have a commutative diagram%
\begin{equation}
\xymatrix{ & & M \times N \ar_{\tau }[ld] \ar^{\beta}[rd] & & \\ & M
\otimes_{S} N \ar[dl]_{{\mathfrak{c}}_{M \otimes_S N}} \ar@{.>}_{\gamma}[rr]
& & G \ar^{{\mathfrak c}_G}[dr] & \\ \mathfrak{c}(M \otimes_{S} N)
\ar@{.>}_{\mathfrak{c}(\gamma)}[rrrr] & & & &\mathfrak{c}(G) }
\label{ten-extend}
\end{equation}%
which suggests that $\mathfrak{c}(-\otimes _{S}-)$ plays the same role.
\end{punto}

The above observations motivates the following connection between the
bifunctors $-\otimes _{S}-$ and $-\boxtimes _{S}-,$ where $\mathbf{CAbMon}$
is the category of cancellative Abelian monoids:

\begin{theorem}
\label{ten-retr}We have an equivalence of functors%
\begin{equation*}
-\boxtimes _{S}-\text{ }\approx \mathfrak{c}(-)\circ (-\otimes _{S}-):%
\mathbb{S}_{S}\times \text{ }_{S}\mathbb{S}\longrightarrow \mathbf{CAbMon}.
\end{equation*}%
In particular, for every right $S$-semimodule $M_{S}$ and every left $S$%
-semimodule $_{S}N,$ we have a natural isomorphism of Abelian monoids%
\begin{equation*}
M\boxtimes _{S}N\simeq \mathfrak{c}(M\otimes _{S}N).
\end{equation*}
\end{theorem}

\begin{Beweis}
Let $M_{S}$ be a right $S$-semimodule, $_{S}N$ a left $S$-semimodule and
consider the Abelian monoids $(M\otimes _{S}N;\tau ),$ $(M\boxtimes _{S}N;%
\widetilde{\tau })$ along with the canonical morphisms of monoids%
\begin{equation*}
\mathfrak{c}_{M\otimes _{S}N}:M\otimes _{S}N\longrightarrow \mathfrak{c}%
(M\otimes _{S}N)\text{ and }\mathfrak{c}(\gamma ):\mathfrak{c}(M\otimes
_{S}N)\longrightarrow M\boxtimes _{S}N.
\end{equation*}%
Since $\widetilde{G}=\mathfrak{c}(M\otimes _{S}N)$ is cancellative and $%
\widetilde{\beta }:=\mathfrak{c}_{M\otimes _{S}N}\circ \tau :M\times
N\longrightarrow \mathfrak{c}(M\otimes N)$ is $S$-balanced, there exists a
\emph{unique} morphism of monoids $\widetilde{\gamma }:M\boxtimes
_{S}N\longrightarrow \mathfrak{c}(M\otimes _{S}N)$ such that $\widetilde{%
\gamma }\circ \widetilde{\tau }=\widetilde{\beta }=\mathfrak{c}_{M\otimes
_{S}N}\circ \tau .$ On other hand, for $G=M\boxtimes _{S}N,$ the map $\beta
:=\widetilde{\tau }:M\times N\longrightarrow M\boxtimes _{S}N$ is $S$%
-balanced and so there exists a \emph{unique} morphism of monoids $\gamma
:M\otimes _{S}N\longrightarrow M\boxtimes _{S}N$ such that $\gamma \circ
\tau =\beta =\widetilde{\tau }.$ Consider the morphisms of monoids%
\begin{eqnarray*}
\varphi &:&M\boxtimes _{S}N\overset{\widetilde{\gamma }}{\longrightarrow }%
\mathfrak{c}(M\otimes _{S}N)\overset{\mathfrak{c}(\gamma )}{\longrightarrow }%
M\boxtimes _{S}N; \\
\theta &:&M\otimes _{S}N\overset{\mathfrak{c}_{M\otimes _{S}N}}{%
\longrightarrow }\mathfrak{c}(M\otimes _{S}N)\overset{\mathfrak{c}(\gamma )}{%
\longrightarrow }M\boxtimes _{S}N\overset{\widetilde{\gamma }}{%
\longrightarrow }\mathfrak{c}(M\otimes _{S}N).
\end{eqnarray*}%
Notice that%
\begin{equation*}
\varphi \circ \widetilde{\tau }=\mathfrak{c}(\gamma )\circ \widetilde{\gamma
}\circ \widetilde{\tau }=\mathfrak{c}(\gamma )\circ \widetilde{\beta }=%
\mathfrak{c}(\gamma )\circ \mathfrak{c}_{M\otimes _{S}N}\circ \tau =\gamma
\circ \tau =\widetilde{\tau }.
\end{equation*}%
Since $\mathrm{id}_{M\boxtimes _{S}N}:M\boxtimes _{S}N\longrightarrow
M\boxtimes _{S}N$ is the \emph{unique} morphism of monoids satisfying this
property, we conclude that $\mathfrak{c}(\gamma )\circ \widetilde{\gamma }=%
\mathrm{id}_{M\boxtimes _{S}N}.$ On the other hand, we have%
\begin{equation*}
\theta \circ \tau =\widetilde{\gamma }\circ \mathfrak{c}(\gamma )\circ
\mathfrak{c}_{M\otimes _{S}N}\circ \tau =\widetilde{\gamma }\circ \gamma
\circ \tau =\widetilde{\gamma }\circ \widetilde{\tau }=\mathfrak{c}%
_{M\otimes _{S}N}\circ \tau .
\end{equation*}%
Although $\tau $ is not an epimorphism (in general), $\tau (M\times N)\ $is
a generating set for $M\otimes _{S}N,$ whence $\theta =\mathfrak{c}%
_{M\otimes _{S}N}$ and so $\widetilde{\gamma }\circ \mathfrak{c}(\gamma
)\circ \mathfrak{c}_{M\otimes _{S}N}=\mathrm{id}_{\mathfrak{c}(M\otimes
_{S}N)}\circ \mathfrak{c}_{M\otimes _{S}N}.$ Since $\mathfrak{c}_{M\otimes
_{S}N}$ is an epimorphism, we conclude that $\widetilde{\gamma }\circ
\mathfrak{c}(\gamma )=\mathrm{id}_{\mathfrak{c}(M\otimes _{S}N)}.$ One can
easily check that this isomorphism is natural in $M_{S}$ and $%
_{S}N.\blacksquare $
\end{Beweis}

\begin{remarks}
Let $S$ and $T$ be semirings.

\begin{enumerate}
\item For every right $S$-semimodule $M$ and left $S$-semimodule $_{S}N$ we
have canonical isomorphisms of Abelian monoids $M\overset{\vartheta _{M}^{r}}%
{\simeq }M\otimes _{S}S$ and $N\overset{\vartheta _{M}^{l}}{\simeq }S\otimes
_{S}N,$ whence%
\begin{equation*}
M\boxtimes _{S}S\simeq \mathfrak{c}(M\otimes _{S}S)\simeq \mathfrak{c}(M)%
\text{ and }S\boxtimes _{S}N\simeq \mathfrak{c}(S\otimes _{S}N)\simeq
\mathfrak{c}(N).
\end{equation*}

\item If $M$ is a right $S$-semimodule, $N$ is an $(S,T)$-bisemimodule and $%
X $ is a left $T$-semimodule, then we have a canonical isomorphism of
Abelian monoids%
\begin{equation*}
(M\otimes _{S}N)\otimes _{T}X\simeq M\otimes _{S}(N\otimes _{T}X).
\end{equation*}
\end{enumerate}
\end{remarks}

\begin{proposition}
\label{adjoint}\emph{(cf. \cite{KN2011}) }Let $M$ be a right $S$-semimodule
and $N$ a left $S$-semimodule.

\begin{enumerate}
\item If $M$ is a $(T,S)$-bisemimodule, then $M\otimes _{S}-:$ $_{S}\mathbb{S%
}\longrightarrow $ $_{T}\mathbb{S}$ is left adjoint to $\mathrm{Hom}%
_{T}(M,-):$ $_{T}\mathbb{S}\longrightarrow $ $_{S}\mathbb{S},$ \emph{i.e.}
for every left $S$-semimodule $X$ and every left $T$-semimodule $Y,$ we have
a canonical isomorphism of Abelian monoids that is natural in $_{S}X$ and $%
_{T}Y:$%
\begin{equation*}
\mathrm{Hom}_{T}(M\otimes _{S}X,Y)\overset{\varsigma ^{l}}{\simeq }\mathrm{%
Hom}_{S}(X,\mathrm{Hom}_{T}(M,Y)).
\end{equation*}

\item If $N$ is an $(S,T)$- bisemimodule, then $-\otimes _{S}N:$ $\mathbb{S}%
_{S}\longrightarrow \mathbb{S}_{T}$ is left adjoint to $\mathrm{Hom}%
_{T}(N,-):\mathbb{S}_{T}\longrightarrow \mathbb{S}_{S},$ \emph{i.e.} for
every right $S$-semimodule $X$ and every right $T$-semimodule $Y$ we have a
canonical isomorphism of Abelian monoids that is natural in $X_{S}$ and $%
Y_{T}:$%
\begin{equation*}
\mathrm{Hom}_{T}(X\otimes _{S}N,Y)\overset{\varsigma ^{r}}{\simeq }\mathrm{%
Hom}_{S}(X,\mathrm{Hom}_{T}(N,Y)).
\end{equation*}
\end{enumerate}
\end{proposition}

As a consequence of Lemma \ref{ten-retr} and Proposition \ref{adjoint} we
recover \cite[Corollary 4.5]{Tak1982c}:

\begin{corollary}
Let $M$ be a right $S$-semimodule and $N$ a left $S$-semimodule.

\begin{enumerate}
\item If $M$ is a $(T,S)$-bisemimodule, $_{S}X$ a left $S$-semimodule and $%
Y\in $ $_{T}\mathbb{CS}$ a $\emph{canecllative}$ left $T$-semimodule, then
we have a canonical isomorphism%
\begin{equation*}
\mathrm{Hom}_{T}(M\boxtimes _{S}X,Y)\simeq \mathrm{Hom}_{S}(X,\mathrm{Hom}%
_{T}(M,Y)).
\end{equation*}

\item If $N$ is an $(S,T)$- bisemimodule, $X$ is a right $S$-semimodule and $%
Y\in \mathbb{CS}_{T}$ a \emph{cancellative }right $T$-semimodule, then we
have a canonical isomorphism%
\begin{equation*}
\mathrm{Hom}_{T}(X\boxtimes _{S}N,Y)\overset{\varsigma ^{r}}{\simeq }\mathrm{%
Hom}_{S}(X,\mathrm{Hom}_{T}(N,Y)).
\end{equation*}
\end{enumerate}
\end{corollary}

\begin{Beweis}
We prove \textquotedblleft 1\textquotedblright . The proof of
\textquotedblleft 2\textquotedblright\ is similar. The required isomorphism
is given by%
\begin{eqnarray*}
\mathrm{Hom}_{T}(M\boxtimes _{S}X,Y) &=&\mathrm{Hom}_{T}(\mathfrak{c}%
(M\otimes _{S}X),Y) \\
&\simeq &\mathrm{Hom}_{T}(M\otimes _{S}X,Y) \\
&\simeq &\mathrm{Hom}_{S}(X,\mathrm{Hom}_{T}(M,Y)).\blacksquare
\end{eqnarray*}
\end{Beweis}

\begin{definition}
A category $\mathfrak{C}$ is said to be (\emph{finitely}) \emph{complete}
iff every functor $F:\mathfrak{D}\longrightarrow \mathfrak{C},$ with $%
\mathfrak{D}$ a small (finite) category, has a limit. Dually, $\mathfrak{C}$
is said to be (\emph{finitely}) \emph{cocomplete} iff every functor $F:%
\mathfrak{D}\longrightarrow \mathfrak{C}$ with $\mathfrak{D}$ a small
(finite) category has a colimit.
\end{definition}

Taking into account the fact that $\mathbb{S}_{S}$ is a variety (in the
sense of Universal Algebra) we have (e.g. \cite[Theorem 21.6.4]{Sch1972}):

\begin{proposition}
\label{cc}The category $\mathbb{S}_{S}$ of right $S$-semimodules is complete
\emph{(}has equalizers and products\emph{)} and cocomplete \emph{(}has
coequalizers and coproducts\emph{)}.
\end{proposition}

\begin{punto}
Let $J$ be a directed set. The \emph{directed {}limit} (\emph{inductive limit%
}, \emph{filtered colimit}) of a \emph{directed system of }$S$\emph{%
-semimodules} $(M_{j},\{f_{jj^{\prime }}:M_{j}\longrightarrow M_{j^{\prime
}}\mid j\leq j^{\prime }\})_{J}$ can be constructed as follows: consider the
disjoint union $\coprod\limits_{j\in J}M_{j}=\bigcup\limits_{j\in
J}(M_{j}\times \{j\}),$ the embeddings $\{\iota _{j}:M_{j}\longrightarrow
\coprod\limits_{j\in J}M_{j}\}_{J}$ and the \emph{congruence relation} on $%
\coprod\limits_{j\in J}M_{j}:$%
\begin{equation}
(x,j)\sim (x^{\prime },j^{\prime })\Longleftrightarrow \text{ }\exists \text{
}j^{\prime \prime }\geq j,j^{\prime }\text{ s.t. }x\in M_{j},\text{ }%
x^{\prime }\in M_{j^{\prime }}\text{ and }f_{jj^{\prime \prime
}}(x)=f_{j^{\prime }j^{\prime \prime }}(x^{\prime }).  \label{dr-eqv}
\end{equation}%
We define%
\begin{equation*}
\lim\limits_{\longrightarrow }M_{j}:=\coprod\limits_{j\in J}M_{j}/\sim \text{
and }\gamma _{j}:M_{j}\longrightarrow \text{ }\lim\limits_{\longrightarrow
}M_{j},\text{ }m\mapsto \lbrack (m,j)].
\end{equation*}%
Notice that $\lim\limits_{\longrightarrow }M_{j}$ is an $S$-semimodule with $%
[(m,j)]s=[(ms,j)]$ for all $s\in S$ and $m\in M_{j}$ and%
\begin{equation*}
\lbrack (m,j)]+[(m^{\prime },j^{\prime })]=[(f_{jj^{\prime \prime
}}(m)+f_{j^{\prime }j^{\prime \prime }}(m^{\prime }),j^{\prime \prime })],%
\text{ where }j^{\prime \prime }\geq j,j^{\prime }.
\end{equation*}
\end{punto}

\begin{punto}
Let $J$ be a directed set. The \emph{inverse {}limit} (\emph{projective limit%
}) of an \emph{inverse system of }$S$\emph{-semimodules} $%
(M_{j},\{f_{jj^{\prime }}:M_{j^{\prime }}\longrightarrow M_{j}\mid j\leq
j^{\prime }\})_{J}$ is given by:%
\begin{equation*}
\lim_{\longleftarrow }M_{j}=\{(m_{j})_{j\in J}\mid \text{ }%
m_{j}=f_{jj^{\prime }}(m_{j^{\prime }})\text{ whenever }j\leq j^{\prime }\}.
\end{equation*}
\end{punto}

The proof of the following important observation is straightforward:

\begin{proposition}
\label{lim-sub}Every $S$-semimodule $M$ is a direct limit of its finitely
generated $S$-subsemimodules.
\end{proposition}

\begin{lemma}
\label{dl-inj}\emph{(\cite{Abu-b})} Let $J$ be a directed set and $%
(X_{j},\{f_{jj^{\prime }}:X_{j}\longrightarrow M_{j^{\prime }}\})_{J},$ $%
(Y_{j},\{g_{jj^{\prime }}:Y_{j}\longrightarrow Y_{j^{\prime }}\})_{J}$ be
directed systems of $S$-semimodules. Let $\{h_{j}:X_{j}\longrightarrow
Y_{j}\}_{J}$ be a class of $S$-linear morphisms satisfying $h_{j^{\prime
}}\circ f_{jj^{\prime }}=g_{jj^{\prime }}\circ h_{j}$ for all $j,j^{\prime
}\in J$ with $j\leq j^{\prime }.$

\begin{enumerate}
\item There exists a unique morphism $h:(\lim\limits_{\longrightarrow
}X_{j},f_{j})\longrightarrow (\lim\limits_{\longrightarrow }Y_{j},g_{j})$
which satisfies $g_{j}\circ h_{j}=h\circ f_{j}.$

\item If $h_{j}$ is injective \emph{(}surjective\emph{)} for every $j\in J,$
then $h$ is injective \emph{(}surjective\emph{)}.

\item If $h_{j}$ is uniform \emph{(}resp. $k$-uniform, $i$-uniform\emph{)}
for every $j\in J,$ then $h$ is uniform \emph{(}resp. $k$-uniform, $i$%
-uniform\emph{)}.
\end{enumerate}
\end{lemma}

\begin{proposition}
\label{ex-d-lim}Let $(L_{j},\{f_{jj^{\prime }}\})_{J},$ $(M_{j},\{g_{jj^{%
\prime }}\})_{J}$ and $(N_{j},\{h_{jj^{\prime }}\})_{J}$ be directed systems
of $S$-semimodules.

\begin{enumerate}
\item If $\{L_{j}\overset{\alpha _{j}}{\longrightarrow }M_{j}\overset{\beta
_{j}}{\longrightarrow }N_{j}\}_{J}$ is a class of exact \emph{(}resp.
semi-exact, proper-exact, quasi-exact\emph{)} sequences of $S$-semimodules,
with $\alpha _{j^{\prime }}\circ f_{jj^{\prime }}=g_{jj^{\prime }}\circ
\alpha _{j}$ and $\beta _{j^{\prime }}\circ g_{jj^{\prime }}=h_{jj^{\prime
}}\circ \beta _{j}$ for all $j\in J,$ then the induced sequence of $S$%
-semimodules $\lim\limits_{\longrightarrow }L_{j}\overset{\alpha }{%
\longrightarrow }\lim\limits_{\longrightarrow }M_{j}\overset{\beta }{%
\longrightarrow }\lim\limits_{\longrightarrow }N_{j}$ is exact \emph{(}resp.
semi-exact, proper-exact, quasi-exact\emph{).}

\item If $\{0\longrightarrow L_{j}\overset{\alpha _{j}}{\longrightarrow }%
M_{j}\overset{\beta _{j}}{\longrightarrow }N_{j}\longrightarrow 0\}_{J}$ is
a class of short exact \emph{(}resp. semi-exact, proper-exact, quasi-exact%
\emph{)} sequences of $S$-semimodules, with $\alpha _{j^{\prime }}\circ
f_{jj^{\prime }}=g_{jj^{\prime }}\circ \alpha _{j}$ and $\beta _{j^{\prime
}}\circ g_{jj^{\prime }}=h_{jj^{\prime }}\circ \beta _{j}$ for all $j\in J,$
then the induced sequence of $S$-semimodules $0\longrightarrow
\lim\limits_{\longrightarrow }L_{j}\overset{\alpha }{\longrightarrow }%
\lim\limits_{\longrightarrow }M_{j}\overset{\beta }{\longrightarrow }%
\lim\limits_{\longrightarrow }N_{j}\longrightarrow 0$ is exact \emph{(}resp.
semi-exact, proper-exact, quasi-exact\emph{). }In particular, $\mathrm{Ker}%
(\beta )\simeq \lim\limits_{\longrightarrow }\mathrm{Ker}(\beta _{j})\ $and $%
\mathrm{Coker}(\alpha )\simeq \lim\limits_{\longrightarrow }\mathrm{Coker}%
(\alpha _{j}).$
\end{enumerate}
\end{proposition}

\begin{lemma}
\label{fg-psi-inj}Let $(M_{j},\{f_{jj^{\prime }}\})_{J}$ be a directed
system of left $S$-semimodules with associated directed system of $S$-linear
maps $f_{j}:M_{j}\longrightarrow \lim\limits_{\longrightarrow }M_{j}$ and
let $X$ be a left $S$-semimodule.

\begin{enumerate}
\item $(\mathrm{Hom}_{S}(X,M_{j}),(X,f_{jj^{\prime }}))_{J}$ is a directed
system of Abelian monoids. Moreover, $(X,f_{j}):\mathrm{Hom}%
_{S}(X,M_{j})\longrightarrow \mathrm{Hom}_{S}(X,\lim\limits_{\longrightarrow
}M_{j})$ is a directed system of morphisms of Abelian monoids and induces a
morphism of Abelian monoids%
\begin{equation}
\psi _{X}=\lim\limits_{\longrightarrow
}(X,f_{j}):\lim\limits_{\longrightarrow }\mathrm{Hom}_{S}(X,M_{j})%
\longrightarrow \mathrm{Hom}_{S}(X,\lim\limits_{\longrightarrow }M_{j}),%
\text{ }[(\alpha _{j},j)]\mapsto \lbrack (f_{j}\circ \alpha _{j},j)]].
\label{X-psi}
\end{equation}

\item If $_{S}X$ is finitely generated, then $\psi _{X}$ is injective.
\end{enumerate}
\end{lemma}

\begin{Beweis}
The first statement is obvious. Assume that $_{S}X$ is finitely generated.
Suppose that $\psi _{X}([(\alpha _{j},j)])=\psi _{X}([(\alpha _{j^{\prime
}}^{\prime },j^{\prime })]),$ \emph{i.e.} $f_{j}\circ \alpha
_{j}=f_{j^{\prime }}\circ \alpha _{j^{\prime }}^{\prime }$ for some $\alpha
_{j}\in \mathrm{Hom}_{S}(X,M_{j}),$ $\alpha _{j}^{\prime }\in \mathrm{Hom}%
_{S}(X,M_{j^{\prime }})$ and $j,j^{\prime }\in J.$ Since $_{S}X$ is finitely
generated, there exists $j^{\prime \prime }\geq j,j^{\prime }$ such that $%
f_{jj^{\prime \prime }}\circ \alpha _{j}=f_{j^{\prime }j^{\prime \prime
}}\circ \alpha _{j^{\prime \prime }},$ \emph{i.e. }$(X,f_{j^{\prime
}j^{\prime \prime }})(\alpha _{j^{\prime }})=(X,f_{jj^{\prime \prime
}})(\alpha _{j}),$ whence $[(\alpha _{j},j)]=[(\alpha _{j^{\prime
}},j^{\prime })].\blacksquare $
\end{Beweis}

\begin{proposition}
\label{adj-lim}\emph{(cf. \cite[Proposition 3.2.2]{Bor1994a})} Let $%
\mathfrak{C},\mathfrak{D}$ be arbitrary categories and $\mathfrak{C}\overset{%
F}{\longrightarrow }\mathfrak{D}\overset{G}{\longrightarrow }\mathfrak{C}$
be functors such that $(F,G)$ is an adjoint pair.

\begin{enumerate}
\item $F$ preserves all colimits which turn out to exist in $\mathfrak{C}.$

\item $G$ preserves all limits which turn out to exist in $\mathfrak{D}.$
\end{enumerate}
\end{proposition}

The following results can be obtained as a direct consequence of
Propositions \ref{adj-lim} and \ref{adjoint}.

\begin{corollary}
\label{ad-l-cor}Let $S,$ $T$ be semirings and $_{T}F_{S}$ a $(T,S)$%
-bisemimodule.

\begin{enumerate}
\item $F\otimes _{S}-:$ $_{S}\mathbb{S}\longrightarrow $ $_{T}\mathbb{S}$
preserves all colimits.

\begin{enumerate}
\item For every family of left $S$-semimodules $\{X_{\lambda }\}_{\Lambda },$
we have a canonical isomorphism of left $T$-semimodules%
\begin{equation*}
F\otimes _{S}\bigoplus\limits_{\lambda \in \Lambda }X_{\lambda }\simeq
\bigoplus\limits_{\lambda \in \Lambda }(F\otimes _{S}X_{\lambda }).
\end{equation*}

\item For any directed system of left $S$-semimodules $(X_{j},\{f_{jj^{%
\prime }}\})_{J},$ we have an isomorphism of left $T$-semimodules%
\begin{equation*}
F\otimes _{S}\lim_{\longrightarrow }X_{j}\simeq \text{ }\lim_{%
\longrightarrow }(F\otimes _{S}X_{j}).
\end{equation*}

\item $F\otimes _{S}-$ preserves coequalizers.

\item $F\otimes _{S}-$ preserves cokernels \emph{(}uniform quotients\emph{)}.
\end{enumerate}

\item $\mathrm{Hom}_{T}(F,-):$ $_{T}\mathbb{S}\longrightarrow $ $_{S}\mathbb{%
S}$ preserves all limits.

\begin{enumerate}
\item For every family of left $T$-semimodules $\{Y_{\lambda }\}_{\Lambda },$
we have a canonical isomorphism of left $S$-semimodules%
\begin{equation*}
\mathrm{Hom}_{T}(F,\prod\limits_{\lambda \in \Lambda }Y_{\lambda })\simeq
\prod\limits_{\lambda \in \Lambda }\mathrm{Hom}_{T}(F,Y_{\lambda }).
\end{equation*}

\item For any inverse system of left $T$-semimodules $(X_{j},\{f_{jj^{\prime
}}\})_{J},$ we have an isomorphism of left $S$-semimodules%
\begin{equation*}
\mathrm{Hom}_{T}(F,\lim_{\longleftarrow }X_{j})\simeq \text{ }%
\lim_{\longleftarrow }\mathrm{Hom}_{T}(F,X_{j}).
\end{equation*}

\item $\mathrm{Hom}_{T}(F,-)$ preserves equalizers;

\item $\mathrm{Hom}_{T}(F,-)$ preserves kernels \emph{(}uniform
subsemimodules\emph{)}.
\end{enumerate}

\item $\mathrm{Hom}_{T}(-,F):$ $_{T}\mathbb{S}\longrightarrow $ $\mathbb{S}%
_{S}$ preserves all limits.

\begin{enumerate}
\item For every family of left $T$-semimodules $\{Y_{\lambda }\}_{\Lambda },$
we have a canonical isomorphism of right $S$-semimodules%
\begin{equation*}
\mathrm{Hom}_{T}(\bigoplus\limits_{\lambda \in \Lambda }Y_{\lambda
},F,)\simeq \prod\limits_{\lambda \in \Lambda }\mathrm{Hom}_{T}(Y_{\lambda
},F).
\end{equation*}

\item For any directed system of left $T$-semimodules $(X_{j},\{f_{jj^{%
\prime }}\})_{J},$ we have an isomorphism of right $S$-semimodules%
\begin{equation*}
\mathrm{Hom}_{T}(\lim_{\longrightarrow }X_{j},F)\simeq \text{ }%
\lim_{\longleftarrow }\mathrm{Hom}_{T}(X_{j},F).
\end{equation*}

\item $\mathrm{Hom}_{T}(-,F)$ converts coequalizers into equalizers;

\item $\mathrm{Hom}_{T}(F,-)$ converts cokernels into kernels \emph{(}%
uniform quotients into uniform subsemimodules\emph{)}.
\end{enumerate}
\end{enumerate}
\end{corollary}

Corollary \ref{ad-l-cor} allows us to improve \cite[Theorem 2.6]{Tak1982a}.

\begin{proposition}
\label{lr-exact}Let $_{T}G_{S}$ an $S$-bisemimodule and consider the functor
$\mathrm{Hom}_{T}(G,-):$ $_{T}\mathbb{S}\longrightarrow $ $_{S}\mathbb{S}.$
Let%
\begin{equation}
0\longrightarrow L\overset{f}{\rightarrow }M\overset{g}{\rightarrow }N
\label{lr}
\end{equation}%
be a sequence of left $T$-semimodules and consider the following sequence of
left $S$-semimodules%
\begin{equation}
0\longrightarrow \mathrm{Hom}_{T}(G,L)\overset{(G,f)}{\rightarrow }\mathrm{%
Hom}_{T}(G,M)\overset{(G,g)}{\longrightarrow }\mathrm{Hom}_{T}(G,N).
\label{GL}
\end{equation}

\begin{enumerate}
\item If $0\longrightarrow L\overset{f}{\rightarrow }M$ is exact and $f$ is
uniform, then $0\longrightarrow \mathrm{Hom}_{T}(G,L)\overset{(G,f)}{%
\rightarrow }\mathrm{Hom}_{T}(G,M)$ is exact and $(G,f)$ is uniform.

\item If \emph{(}\ref{lr}\emph{)\ }is semi-exact and $f$ is uniform, then
\emph{(}\ref{GL}\emph{)\ }is semi-exact \emph{(}proper exact\emph{)} and $%
(G,f)$ is uniform.

\item If \emph{(}\ref{lr}\emph{)\ }is exact and $\mathrm{Hom}_{T}(G,-)$
preserves $k$-uniform morphisms, then \emph{(}\ref{GL}\emph{)\ }is exact.
\end{enumerate}
\end{proposition}

\begin{Beweis}
\begin{enumerate}
\item The following implications are obvious: $0\longrightarrow L\overset{f}{%
\rightarrow }M$ is exact $\Longrightarrow $ $f$ is injective $%
\Longrightarrow $ $(G,f)$ is injective $\Longrightarrow 0\longrightarrow
\mathrm{Hom}_{T}(G,L)\overset{(G,f)}{\rightarrow }\mathrm{Hom}_{T}(G,M)$ is
exact. Assume that $f$ is uniform and consider the exact sequence of $S$%
-semimodules%
\begin{equation*}
0\longrightarrow L\overset{f}{\longrightarrow }M\overset{\pi _{L}}{%
\longrightarrow }M/L\longrightarrow 0.
\end{equation*}%
Notice that $L=\mathrm{Ker}(\pi _{L})$ by Lemma \ref{exact} (5). By
Corollary \ref{ad-l-cor}, $\mathrm{Hom}_{T}(G,-)$ preserves kernels and so $%
(G,f)=\mathrm{ker}(G,\pi _{L})$ whence uniform.

\item Apply Lemma \ref{exact} (3): The semi-exactness of (\emph{\ref{lr}})
and the uniformity of $f$ are equivalent to $L\simeq \mathrm{Ker}(g).$ Since
$\mathrm{Hom}_{T}(G,-)$ preserves kernels, we deduce that $\mathrm{Hom}%
_{T}(G,L)=\mathrm{Ker}((G,g))$ which is equivalent to the semi-exactness of (%
\emph{\ref{GL}}) and the uniformity of $(G,f).$ Notice that $(G,f)(\mathrm{%
Hom}_{T}(G,L))=\overline{(G,f)(\mathrm{Hom}_{T}(G,L))}=\mathrm{Ker}(G,g),$
\emph{i.e. }(\emph{\ref{GL}}) is proper exact.

\item The statement follows directly from \textquotedblleft
2\textquotedblright\ and the assumption on $\mathrm{Hom}_{T}(G,-).%
\blacksquare $
\end{enumerate}
\end{Beweis}

\begin{proposition}
\label{ll-exact}Let $_{T}G_{S}$ be a $(T,S)$-bisemimodule and consider the
functor $\mathrm{Hom}_{T}(-,G):$ $_{T}\mathbb{S}\longrightarrow \mathbb{S}%
_{S}.$ Let%
\begin{equation}
L\overset{f}{\rightarrow }M\overset{g}{\rightarrow }N\longrightarrow 0
\label{L-l}
\end{equation}%
be a sequence of left $T$-semimodules and consider the sequence of right $S$%
-semimodules%
\begin{equation}
0\longrightarrow \mathrm{Hom}_{T}(N,G)\overset{(g,G)}{\rightarrow }\mathrm{%
Hom}_{T}(M,G)\overset{(f,G)}{\longrightarrow }\mathrm{Hom}_{T}(L,G).
\label{LG}
\end{equation}

\begin{enumerate}
\item If $M\overset{g}{\rightarrow }N\longrightarrow 0$ is exact and $g$ is
uniform, then $0\longrightarrow \mathrm{Hom}_{T}(N,G)\overset{(g,G)}{%
\rightarrow }\mathrm{Hom}_{T}(M,G)$ is exact and $(g,G)$ is uniform.

\item If \emph{(\ref{L-l}) }is semi-exact and $g$ is uniform, then \emph{(%
\ref{LG}) }is semi-exact \emph{(}proper-exact\emph{)} and $(g,G)$ is uniform.

\item If \emph{(\ref{L-l}) }is exact and $\mathrm{Hom}_{T}(-,G)$ converts $i$%
-uniform morphisms into $k$-uniform ones, then \emph{(\ref{LG}) }is exact.
\end{enumerate}
\end{proposition}

\begin{Beweis}
\begin{enumerate}
\item The following implications are clear: $M\overset{g}{\rightarrow }%
N\longrightarrow 0$ is exact $\Longrightarrow $ $g$ is surjective $%
\Longrightarrow $ $(g,G)$ is injective $\Longrightarrow $ $0\longrightarrow
\mathrm{Hom}_{T}(N,G)\overset{(g,G)}{\rightarrow }\mathrm{Hom}_{T}(M,G)$ is
exact. Assume that $g$ is uniform and consider the exact sequence of $S$%
-semimodules%
\begin{equation*}
0\longrightarrow \mathrm{Ker}(g)\overset{\iota }{\longrightarrow }M\overset{g%
}{\longrightarrow }N\longrightarrow 0.
\end{equation*}%
Notice that $N\simeq \mathrm{Coker}(\iota ).$ By Corollary \ref{ad-l-cor}, $%
\mathrm{Hom}_{T}(-,G)$ converts cokernels into kernels, we conclude that $%
(g,G)=\mathrm{ker}((f,G))$ whence uniform.

\item Apply Lemma \ref{exact}: $L\overset{f}{\rightarrow }M\overset{g}{%
\rightarrow }N\longrightarrow 0$ is semi-exact and $f$ is uniform $%
\Longleftrightarrow $ $M\simeq \mathrm{Coker}(f).$ Since the contravariant
functor $\mathrm{Hom}_{T}(-,G)$ converts cokernels into kernels, it follows
that $\mathrm{Hom}_{T}(N,G)=\mathrm{Ker}((f,G))$ which is in turn equivalent
to (\emph{\ref{LG}}) being semi-exact and $(g,G)$ being uniform. Notice that
$(g,G)(\mathrm{Hom}_{S}(N,G))=\overline{(g,G)(\mathrm{Hom}_{S}(N,G))}=%
\mathrm{Ker}((f,G)),$ \emph{i.e.} (\emph{\ref{LG}}) is proper-exact.

\item This follows immediately from \textquotedblleft 2\textquotedblright\
and the assumption on $\mathrm{Hom}_{T}(-,G).\blacksquare $
\end{enumerate}
\end{Beweis}

\begin{proposition}
\label{r-exact}Let $_{T}G_{S}$ be a $(T,S)$-bisemimodule and consider the
functor $G\otimes _{S}-:$ $_{S}\mathbb{S}\longrightarrow $ $_{T}\mathbb{S}.$
Let%
\begin{equation}
L\overset{f}{\rightarrow }M\overset{g}{\rightarrow }N\rightarrow 0
\label{L-r}
\end{equation}%
be a sequence of left $S$-semimodules and consider the sequence of left $T$%
-semimodules%
\begin{equation}
G\otimes _{S}L\overset{\mathrm{id}_{G}\otimes _{S}f}{\rightarrow }G\otimes
_{S}M\overset{\mathrm{id}_{G}\otimes _{S}g}{\longrightarrow }G\otimes
_{S}N\rightarrow 0  \label{FL-r}
\end{equation}

\begin{enumerate}
\item If $M\overset{g}{\rightarrow }N\rightarrow 0$ is exact and $g$ is
uniform, then $G\otimes _{S}M\overset{\mathrm{id}_{G}\otimes _{S}g}{%
\longrightarrow }G\otimes _{S}N\rightarrow 0$ is exact and $\mathrm{id}%
_{G}\otimes _{S}g$ is uniform.
\end{enumerate}
\end{proposition}

\begin{proposition}
If \emph{(\ref{L-r}) }is semi-exact and $g$ is uniform, then \emph{(\ref{L-r}%
) }is semi-exact and $\mathrm{id}_{G}\otimes _{S}g$ is uniform.
\end{proposition}

\begin{proposition}
If \emph{(\ref{L-r}) }is exact and $G\otimes _{S}-$ preserves $i$-uniform
morphisms, then \emph{(\ref{L-r}) }is exact.
\end{proposition}

\begin{Beweis}
The following implications are obvious: $M\overset{g}{\rightarrow }%
N\rightarrow 0$ is exact $\Longrightarrow $ $g$ is surjective $%
\Longrightarrow $ $\mathrm{id}_{G}\otimes _{S}g$ is surjective $%
\Longrightarrow $ $G\otimes _{S}M\overset{\mathrm{id}_{G}\otimes _{S}g}{%
\longrightarrow }G\otimes _{S}N\rightarrow 0$ is exact. Assume that $g$ is
uniform and consider the exact sequence of $S$-semimodules%
\begin{equation*}
0\longrightarrow \mathrm{Ker}(g)\overset{\iota }{\longrightarrow }M\overset{g%
}{\longrightarrow }N\longrightarrow 0.
\end{equation*}%
Then $N\simeq \mathrm{Coker}(\iota ).$ By Corollary \ref{ad-l-cor}, $%
G\otimes _{S}-$ preserves cokernels and so $\mathrm{id}_{G}\otimes _{S}g=%
\mathrm{coker}(\mathrm{id}_{G}\otimes _{S}\iota )$ whence uniform.
\end{Beweis}

\begin{Beweis}
Apply Lemma \ref{exact}: The assumptions on (\ref{L-r}) are equivalent to $N=%
\mathrm{Co}$\textrm{$k$}$\mathrm{er}(f)$ by Lemma \ref{exact}. Since $%
G\otimes _{S}-$ preserves cokernels, we conclude that $G\otimes _{S}N=%
\mathrm{Coker}(\mathrm{id}_{G}\otimes _{S}f),$ \emph{i.e.} (\ref{FL-r}) is
semi-exact and $\mathrm{id}_{G}\otimes _{S}g$ is uniform.
\end{Beweis}

\begin{Beweis}
This follows directly form \textquotedblleft 2\textquotedblright\ and the
assumption on $G\otimes _{S}-.\blacksquare $
\end{Beweis}

We say that an $S$-semimodule $P$ is \emph{projective} iff for every
surjective morphism of $S$-semimodules $M\overset{g}{\longrightarrow }%
N\longrightarrow 0,$ the induced morphism of Abelian monoids $\mathrm{Hom}%
_{S}(P,M)\overset{(P,g)}{\longrightarrow }\mathrm{Hom}_{S}(P,N)%
\longrightarrow 0$ is surjective. It is well-known that $_{S}P$ is
projective if and only if $P$ is a retract of a free $S$-semimodule (\emph{%
e.g.} \cite[Theorem 1.9]{Tak1983}, \cite[Proposition 17.16]{Go19l99a}).

The proof of the following lemma is straightforward:

\begin{lemma}
\label{u-sum}

\begin{enumerate}
\item Let $\{f_{\lambda }:L_{\lambda }\longrightarrow M_{\lambda
}\}_{\Lambda }$ be a family of left $S$-semimodule morphisms and consider
the induced $S$-linear map $f:\bigoplus\limits_{\lambda \in \Lambda
}L_{\lambda }\longrightarrow \bigoplus\limits_{\lambda \in \Lambda
}M_{\lambda }.$ Then $f$ is uniform \emph{(}resp. $k$-uniform, $i$-uniform%
\emph{)} if and only if $f_{\lambda }$ is uniform \emph{(}resp. $k$-uniform,
$i$-uniform\emph{)} for every $\lambda \in \Lambda .$ In particular, $%
\bigoplus\limits_{\lambda \in \Lambda }L_{\lambda }\leq
_{S}^{u}\bigoplus\limits_{\lambda \in \Lambda }M_{\lambda }$ if and only if $%
L_{\lambda }\leq _{S}^{u}M_{\lambda }$ for every $\lambda \in \Lambda .$

\item A morphism $\varphi :L\longrightarrow M$ of left $S$-semimodules is
uniform \emph{(}resp. $k$-uniform, $i$-uniform\emph{)} if and only if $%
\mathrm{id}_{F}\otimes _{S}\varphi :F\otimes _{S}L\longrightarrow F\otimes
_{S}M$ is uniform \emph{(}resp. $k$-uniform, $i$-uniform\emph{)} for every
non-zero free right $S$-semimodule $F\neq 0.$

\item If $P_{S}$ is projective and $\varphi :L\longrightarrow M$ is a
uniform \emph{(}resp. $k$-uniform, $i$-uniform\emph{)} morphism of left $S$%
-semimodules, then $\mathrm{id}_{F}\otimes _{S}\varphi :P\otimes
_{S}L\longrightarrow P\otimes _{S}M$ is uniform \emph{(}resp. $k$-uniform, $%
i $-uniform\emph{)}.
\end{enumerate}
\end{lemma}

It is well-known, that for every (finitely generated) $S$-semimodule $X,$
there is a free $S$-semimodule $S^{(J)},$ for some (finite) index set $J,$
and a surjective $S$-linear map $S^{(J)}\longrightarrow X\longrightarrow 0.$

\begin{definition}
We call a left $S$-semimodule $X:$

\emph{uniformly finitely generated} iff there exists a uniform surjective $S$%
-linear map $S^{n}\longrightarrow X\overset{g}{\longrightarrow }0$ for some $%
n\in \mathbb{N};$

\emph{uniformly finitely presented} iff $_{S}X$ is uniformly finitely
generated and for any exact sequence of $S$-semimodules%
\begin{equation*}
0\longrightarrow K\overset{f}{\longrightarrow }S^{n}\overset{g}{%
\longrightarrow }X\longrightarrow 0,
\end{equation*}%
the $S$-semimodule $K$ is finitely generated.
\end{definition}

\begin{remark}
Takahashi \cite{Tak1983} defined an $S$-semimodule $X$ to be \emph{normal}
iff there exists a projective $S$-semimodule $P$ and a uniform surjective $S$%
-linear map $P\overset{\varepsilon }{\longrightarrow }X\longrightarrow 0$
(called a \emph{projective presentation of }$X$). Indeed, every uniformly
finitely generated $S$-semimodule is normal.
\end{remark}

\begin{proposition}
\label{fp-m-n}If $_{S}X$ is uniformly finitely presented, then there exist $%
m,n\in \mathbb{N}$ and an exact sequence of $S$-semimodules%
\begin{equation*}
S^{m}\overset{\tilde{f}}{\longrightarrow }S^{n}\overset{\tilde{g}}{%
\longrightarrow }X\longrightarrow 0.
\end{equation*}
\end{proposition}

\begin{Beweis}
Since $_{S}X$ is uniformly finitely generated, there exists a uniform
surjective $S$-linear map $\widetilde{g}:S^{n}\longrightarrow X.$ Let $K=%
\mathrm{Ker}(g)$ and consider the exact sequence of left $S$-semimodules%
\begin{equation*}
0\longrightarrow K\overset{\mathrm{ker}(g)}{\longrightarrow }S^{n}\overset{%
\widetilde{g}}{\longrightarrow }X\longrightarrow 0.
\end{equation*}%
By assumption, $_{S}K$ is finitely generated and so there exists a
surjective $S$-linear map $\pi :S^{m}\longrightarrow K$ for some $m\in
\mathbb{N}.$ Notice that $\widetilde{f}:=\mathrm{ker}(g)\circ \pi $ is $i$%
-uniform by \cite[Lemma 3.8 \textquotedblleft 1-c\textquotedblright ]{Abu}
and $\widetilde{g}$ is uniform by assumption. Indeed, $\widetilde{f}(S^{m})=%
\mathrm{Ker}(\widetilde{g})$ and so the following sequence is exact%
\begin{equation*}
S^{m}\overset{\tilde{f}}{\longrightarrow }S^{n}\overset{\tilde{g}}{%
\longrightarrow }X\longrightarrow 0.\blacksquare
\end{equation*}
\end{Beweis}

\begin{definition}
(\cite{Abu-b})\ We say that a right $S$-semimodule $Q$ is (\emph{uniformly})
$\mathcal{M}$-\emph{injective}, where $\mathcal{M}$ is a class of right $S$%
-semimodules, iff for every (uniform) injective morphism $0\longrightarrow L%
\overset{f}{\longrightarrow }M$ with $M\in \mathcal{M},$ the induced
morphism of Abelian monoids $\mathrm{Hom}_{S}(M,Q)\overset{(f,Q)}{%
\longrightarrow }\mathrm{Hom}_{S}(L,Q)$ is surjective (and uniform). If $%
_{S}Q$ is (uniformly) $M$-injective for every $M\in $ $_{S}\mathbb{S},$ then
we say that $_{S}Q$ is (\emph{uniformly}) \emph{injective}. In fact, $_{S}Q$
is uniformly injective if and only if $\mathrm{Hom}_{S}(-,Q)$ preserves
exact sequences.
\end{definition}

\section{Flat Semimodules}

\qquad As before, $S$ is a semiring with $1_{S}\neq 0_{S}.$ If $M$ is a left
$S$-semimodule, then we write $U\leq _{S}^{u}M$ to indicate that $U$ is a
\emph{uniform} (\emph{subtractive}) $S$-subsemimodule of $M$ (\emph{i.e. }%
the embedding map $U\overset{\iota }{\hookrightarrow }M$ is uniform).

\qquad The following definition applies to any variety in the sense of
Universal Algebra (\emph{e.g.} \cite{BR2004}):

\begin{definition}
We say that a right $S$-semimodule $F$ is \emph{flat} iff $%
F=\lim\limits_{\longrightarrow }F_{i},$ a directed limit (filtered colimit)
of finitely presented projective right $S$-semimodules.
\end{definition}

\begin{lemma}
\emph{(cf. \cite{Kat2004a})} The following are equivalent for a right $S$%
-semimodule $F_{S}:$

\begin{enumerate}
\item $F\otimes _{S}-$ is left exact \emph{(i.e.} preserves finite limits%
\emph{)};

\item $F\otimes _{S}-$ preserves pullbacks and equalizers;

\item $F_{S}$ is pullback-flat, i.e. $F\otimes _{S}-$ preserves pullbacks;

\item $F_{S}$ is $L$-flat, i.e. $F\simeq \lim\limits_{\longrightarrow
}F_{\lambda },$ a filtered \emph{(}directed\emph{)} colimit of finitely
generated free $S$-semimodules;

\item $F_{S}$ is flat.
\end{enumerate}
\end{lemma}

Although the above definition is quite natural, a notion of flatness w.r.t.
to a family of semimodules is important. This motivates introducing the
following notion.

\begin{definition}
Let $F$ be a right $S$-semimodule and $\mathcal{M}$ a class of left $S$%
-semimodules. We say that $F$ is \emph{uniformly flat w.r.t.} $\mathcal{M}$
(or \emph{uniformly }$\mathcal{M}$\emph{-flat}) iff for every exact sequence
of left $S$-semimodules%
\begin{equation*}
0\longrightarrow L\overset{f}{\longrightarrow }M\overset{g}{\longrightarrow }%
N\longrightarrow 0,
\end{equation*}%
with $M\in \mathcal{M},$ the following sequence of Abelian monoids is exact%
\begin{equation}
0\longrightarrow F\otimes _{S}L\overset{\mathrm{id}_{F}\otimes _{S}f}{%
\longrightarrow }F\otimes _{S}M\overset{\mathrm{id}_{F}\otimes _{S}g}{%
\longrightarrow }F\otimes _{S}N\longrightarrow 0.  \label{ten-seq}
\end{equation}%
If $F_{S}$ is uniformly $M$-flat for every left $S$-semimodule $_{S}M,$ then
we say that $F$ is \emph{uniformly flat}.
\end{definition}

\begin{theorem}
\label{u-flat}Let $F$ be a right $S$-semimodule.

\begin{enumerate}
\item Let $_{S}M$ be a left $S$-semimodule. Then $F_{S}$ is uniformly $M$%
-flat if and only if for every $U\leq _{S}^{u}M$ we have $F\otimes _{S}U\leq
_{S}^{u}F\otimes _{S}M.$

\item $F_{S}$ is uniformly flat if and only if $F_{S}\otimes -$ preserves
uniform subsemimodule.
\end{enumerate}
\end{theorem}

\begin{Beweis}
We need only to prove \textquotedblleft 1\textquotedblright .

$(\Longrightarrow )$ Assume that $F_{S}$ is uniformly $M$-flat. Let $U\leq
_{S}^{u}M$ and consider the exact sequence of $S$-semimodules $%
0\longrightarrow U\overset{\iota }{\longrightarrow }M\overset{\pi }{%
\longrightarrow }M/U\longrightarrow 0,$ where $\iota $ is the canonical
embedding and $\pi $ is the canonical uniform surjection. By assumption, the
sequence $0\longrightarrow F\otimes _{S}U\overset{\mathrm{id}_{F}\otimes
_{S}\iota }{\longrightarrow }F\otimes _{S}M\overset{\mathrm{id}_{F}\otimes
_{S}\pi }{\longrightarrow }F\otimes _{S}M/U\longrightarrow 0$ is exact; in
particular, $F\otimes _{S}U\leq ^{u}F\otimes _{S}M$ is a uniform submonoid.

$(\Longleftarrow )$ Let $0\longrightarrow L\overset{f}{\longrightarrow }M%
\overset{g}{\longrightarrow }N\longrightarrow 0$ be an exact sequence of
left $S$-semimodules, \emph{i.e.} $L\simeq \mathrm{Ker}(g)$ and $N\simeq
\mathrm{Coker}(f).$ By Proposition \ref{r-exact} \textquotedblleft
2\textquotedblright , the sequence $F\otimes _{S}L\overset{\mathrm{id}%
_{F}\otimes _{S}f}{\longrightarrow }F\otimes _{S}M\overset{\mathrm{id}%
_{F}\otimes _{S}g}{\longrightarrow }F\otimes _{S}N\longrightarrow 0$ is
proper exact and $\mathrm{id}_{F}\otimes _{S}g$ is uniform. By assumption, $%
\mathrm{id}_{F}\otimes _{S}f$ is injective and uniform, whence (\ref{ten-seq}%
) is exact.$\blacksquare $
\end{Beweis}

\begin{corollary}
\label{retr}

\begin{enumerate}
\item Let $M$ be a left $S$-semimodule. Any retract of a uniformly $M$-flat
right $S$-semimodule is uniformly $M$-flat.

\item Any retract of a uniformly flat right $S$-semimodule is uniformly flat.
\end{enumerate}
\end{corollary}

\begin{Beweis}
We need only to prove \textquotedblleft 1\textquotedblright . Let $M$ be a
left $S$-semimodule and $U\leq _{S}^{u}M.$ Let $F_{S}$ be a uniformly $M$%
-flat right $S$-semimodule and $\widetilde{F}$ a retract of $F.$ Then there
exist $S$-linear maps $\widetilde{F}\overset{\psi }{\longrightarrow }F%
\overset{\theta }{\longrightarrow }\widetilde{F}$ such that $\theta \circ
\psi =\mathrm{id}_{\widetilde{F}}.$ Consider the commutative diagram%
\begin{equation*}
\xymatrix{\tilde{F} \otimes_S U \ar[rr]^{{\rm id}_{\tilde{F}}\otimes_S
\iota_U} \ar[d]_{\psi \otimes_S {\rm id}_U} & & \tilde{F} \otimes_S M
\ar[d]^{\psi \otimes_S {\rm id}_M}\\ F \otimes_S U \ar[rr]^{{\rm
id}_{F}\otimes_S \iota_U} \ar[d]_{\theta \otimes_S {\rm id}_U} & & F
\otimes_S M \ar[d]^{\theta \otimes_S {\rm id}_M} \\ \tilde{F} \otimes_S U
\ar[rr]^{{\rm id}_{\tilde{F}}\otimes_S \iota_U} & & \tilde{F} \otimes_S M}
\end{equation*}%
Indeed, $(\theta \otimes _{S}\mathrm{id}_{U})\circ (\psi \otimes _{S}\mathrm{%
id}_{U})=\mathrm{id}_{\widetilde{F}\otimes _{S}U}$ and $(\theta \otimes _{S}%
\mathrm{id}_{M})\circ (\psi \otimes _{S}\mathrm{id}_{M})=\mathrm{id}_{%
\widetilde{F}\otimes _{S}M},$ \emph{i.e.} $\widetilde{F}\otimes _{S}U$ is a
retract of $F\otimes _{S}U$ and $\widetilde{F}\otimes _{S}M$ is a retract of
$F\otimes _{S}M.$ Since $F_{S}$ is flat, $\mathrm{id}_{F}\otimes _{S}\iota
_{U}:F\otimes _{S}U\longrightarrow F\otimes _{S}M$ is injective and uniform.
It follows that $\mathrm{id}_{\widetilde{F}}\otimes _{S}\iota _{U}$ is
injective and indeed uniform by Lemma \ref{comm} \textquotedblleft
1\textquotedblright , \emph{i.e.} $\widetilde{F}\otimes _{S}U\leq _{S}^{u}%
\widetilde{F}\otimes _{S}M.$ Consequently, $\widetilde{F}$ is uniformly $M$%
-flat.$\blacksquare $
\end{Beweis}

\begin{proposition}
\label{sum-flat}Let $\{F_{\lambda }\}_{\Lambda }$ be a family of right $S$%
-semimodules.

\begin{enumerate}
\item Let $M$ be a left $S$-semimodule. Then $\bigoplus\limits_{\lambda \in
\Lambda }F_{\lambda }$ is uniformly $M$-flat if and only if $F_{\lambda }$
is uniformly $M$-flat for every $\lambda \in \Lambda .$

\item $\bigoplus\limits_{\lambda \in \Lambda }F_{\lambda }$ is uniformly
flat if and only if $F_{\lambda }$ is uniformly flat for every $\lambda \in
\Lambda .$
\end{enumerate}
\end{proposition}

\begin{Beweis}
We need only to prove \textquotedblleft 1\textquotedblright . Let $F:=$ $%
\bigoplus\limits_{\lambda \in \Lambda }F_{\lambda }$ and consider the
projections $\pi _{\lambda }:F\longrightarrow F_{\lambda },$ $(f_{\lambda
})_{\Lambda }\mapsto f_{\lambda }$ for $\lambda \in \Lambda .$ Let $U\leq
_{S}^{u}M$ be a uniform $S$-subsemimodule. Assume that $F_{\lambda }$ is $M$%
-flat for every $\lambda \in \Lambda .$ Then $F_{\lambda }\otimes _{S}U\leq
_{S}^{u}F_{\lambda }\otimes _{S}M$ for every $\lambda \in \Lambda ,$ whence $%
\bigoplus\limits_{\lambda \in \Lambda }(F_{\lambda }\otimes _{S}U)\leq
_{S}^{u}\bigoplus\limits_{\lambda \in \Lambda }(F_{\lambda }\otimes _{S}M)$
by Lemma \ref{u-sum}. Since $\bigoplus\limits_{\lambda \in \Lambda
}(F_{\lambda }\otimes _{S}U)\simeq \bigoplus\limits_{\lambda \in \Lambda
}F_{\lambda }\otimes _{S}U$ and $\bigoplus\limits_{\lambda \in \Lambda
}(F_{\lambda }\otimes _{S}M)\simeq \bigoplus\limits_{\lambda \in \Lambda
}F_{\lambda }\otimes _{S}M,$ we conclude that $\bigoplus\limits_{\lambda \in
\Lambda }F_{\lambda }\otimes _{S}U\leq _{S}^{u}\bigoplus\limits_{\lambda \in
\Lambda }F_{\lambda }\otimes _{S}M.$ It follows that $\bigoplus\limits_{%
\lambda \in \Lambda }F_{\lambda }\otimes $ is uniformly $M$-flat. On the
other hand, assume that $\bigoplus\limits_{\lambda \in \Lambda }F_{\lambda }$
is uniformly $M$-flat. Every $F_{\lambda },$ $\lambda \in \Lambda ,$ is a
retract of\emph{\ }$\bigoplus\limits_{\lambda \in \Lambda }F_{\lambda }$
whence uniformly $M$-flat by Corollary \ref{retr}.$\blacksquare $
\end{Beweis}

\begin{lemma}
\label{fp-iso}Let $_{S}X$ be a left $S$-semimodule, $_{S}Y_{T}$ an $(S,T)$%
-bisemimodule, $_{T}Z$ a uniformly flat left $T$-module and consider the
following map of Abelian monoids%
\begin{equation*}
\nu _{X,Y,Z}:\mathrm{Hom}_{S}(X,Y)\otimes _{T}Z\longrightarrow \mathrm{Hom}%
_{S}(X,Y\otimes _{T}Z),\text{ }f\otimes _{T}z\mapsto f(-)\otimes _{S}z].
\end{equation*}

\begin{enumerate}
\item If $_{S}X$ is uniformly finitely generated, then $\nu _{X,Y,Z}$ is
injective and uniform.

\item If $_{S}X$ is uniformly finitely presented, then $\nu _{X,Y,Z}$ is an
isomorphism.
\end{enumerate}
\end{lemma}

\begin{Beweis}
\begin{enumerate}
\item Since $_{S}X$ is uniformly finitely generated, there exists a uniform
surjective $S$-linear map%
\begin{equation*}
S^{n}\overset{\widetilde{g}}{\longrightarrow }X\longrightarrow 0.
\end{equation*}%
By Proposition \ref{ll-exact}, $\mathrm{Hom}_{S}(X,Y)\leq _{S}^{u}\mathrm{Hom%
}_{S}(S^{n},Y),$ whence $\mathrm{Hom}_{S}(X,Y)\otimes _{T}Z\overset{(%
\widetilde{g},Y)\otimes _{T}\mathrm{id}_{Z}}{\hookrightarrow }\mathrm{Hom}%
_{S}(S^{n},Y)\otimes _{T}Z$ since $_{T}Z$ is uniformly flat and we have a
commutative diagram%
\begin{equation*}
\xymatrix{0 \ar[r] & \mathrm{Hom}_{S}(X,Y)\otimes_{T} Z
\ar[rr]^{({\tilde{g}},Y)\otimes_T {\rm id}_Z }\ar[d]_{\nu _{X,Y,Z}} & &
\mathrm{Hom}_{S}(S^n,Y)\otimes _{T}Z \ar[d]_{\nu_{S^n,Y,Z}} \\ 0 \ar[r] &
\mathrm{Hom}_{S}(X,Y\otimes _{T}Z) \ar[rr]_{({\tilde{g}}, Y \otimes_T Z) } &
& \mathrm{Hom}_{S}(S^n,Y \otimes _{T}Z) }
\end{equation*}%
Notice that $\nu _{S^{n},Y,Z}$ is an isomorphism, whence $\nu _{X,Y,Z}$ is
injective. Moreover, it follows by \cite[Lemma 3.8 (1)]{Abu} that $\nu
_{X,Y,Z}$ is uniform.

\item Since $_{S}X$ is finitely presented, there exists by Proposition \ref%
{fp-m-n} an exact sequence of $S$-semimodules $S^{m}\overset{\tilde{f}}{%
\longrightarrow }S^{n}\overset{\tilde{g}}{\longrightarrow }X\longrightarrow
0 $ for some $m,n\in \mathbb{N}.$ By Proposition \ref{ll-exact} and the
uniform flatness of $_{T}Z$ we obtain the following commutative diagram with
proper-exact rows%
\begin{equation*}
\xymatrix{0 \ar[r] & \mathrm{Hom}_{S}(X,Y)\otimes_{T} Z
\ar[rr]^{({\tilde{g}},Y)\otimes_T {\rm id}_Z }\ar[d]_{\nu _{X,Y,Z}} & &
\mathrm{Hom}_{S}(S^n,Y)\otimes _{T}Z \ar[rr]^{({\tilde{f}},Y)\otimes_T {\rm
id}_Z } \ar[d]_{\nu_{S^n,Y,Z}} & & \mathrm{Hom}_{S}(S^m,Y)\otimes_{T} Z
\ar[d]_{\nu_{S^m,Y,Z}} \\ 0 \ar[r] & \mathrm{Hom}_{S}(X,Y\otimes _{T}Z)
\ar[rr]_{({\tilde{g}}, Y \otimes_T Z) } & & \mathrm{Hom}_{S}(S^n,Y \otimes
_{T}Z) \ar[rr]_{({\tilde{f}}, Y \otimes_T Z)} & & \mathrm{Hom}_{S}(S^m,Y
\otimes _{T}Z)}
\end{equation*}%
Notice that $\nu _{S^{m},Y,Z}$ and $\nu _{S^{n},Y,Z}$ are isomorphisms and
so it follows by Lemma \ref{diagram} \textquotedblleft 3\textquotedblright\
that $\nu _{X,Y,Z}$ is surjective. Notice that $\nu _{X,Y,Z}$ is injective
by \textquotedblleft 1\textquotedblright\ whence an isomorphism.$%
\blacksquare $
\end{enumerate}
\end{Beweis}

\qquad Applying Lemma \ref{fp-iso} to $S=T$ and $Y=S,$ considered as a
bisemimodule over itself in the canonical way, we obtain with $X^{\ast }=%
\mathrm{Hom}_{S}(X,S):$

\begin{proposition}
\label{fp-iso-R}Let $_{S}X$ be a uniformly finitely presented $S$%
-semimodule, $_{S}Z$ a uniformly flat left $S$-semimodule and consider the
following morphism of Abelian monoids%
\begin{equation*}
\nu _{X,Z}:X^{\ast }\otimes _{S}Z\longrightarrow \mathrm{Hom}_{S}(X,Z),\text{
}f\otimes _{S}z\mapsto f(-)z].
\end{equation*}%
If $_{S}X$ is uniformly finitely generated \emph{(}uniformly finitely
presented\emph{)}, then $\nu _{X,Z}$ is injective and uniform \emph{(}an
isomorphism\emph{)}.
\end{proposition}

\begin{definition}
Let $M$ be a left $S$-semimodule. We say that a right $S$-semimodule $F_{S}$
is $M$-\emph{mono-flat} \cite{Kat2004a} (or $M$-$k$\emph{-flat }\cite%
{Alt2004}) iff $F\otimes _{S}L\leq _{S}F\otimes _{S}M$ for every $S$%
-subsemimodule $L\leq _{S}M.$ If $F_{S}$ is $M$-mono-flat for every left $S$%
-semimodule $M,$ then we call $F_{S}$ \emph{mono-flat} (or $k$\emph{-flat}).
\end{definition}

\begin{notation}
For every left $S$-semimodule $M,$ we set%
\begin{equation*}
\mathcal{I}_{S}(M):=\{G\in \mathbb{S}_{S}\mid G\otimes _{S}U\overset{\mathrm{%
id}_{G}\otimes _{S}\iota }{\longrightarrow }G\otimes _{S}M\text{ is }i\text{%
-uniform }\forall \text{ }U\leq _{S}^{u}M\};
\end{equation*}
\end{notation}

\begin{remark}
Let $M$ be a left $S$-semimodule. If $F_{S}\in \mathcal{I}_{S}(M)$ and $M$%
-mono-flat, then $F$ is uniformly $M$-flat.
\end{remark}

The following result is straightforward (cf. \cite[Proposition 4.1]{Alt2004}%
):

\begin{proposition}
\label{fg-M-flat}Let $M$ be a left $S$-semimodule and $F_{S}\in \mathcal{I}%
_{S}(M).$ Then $F$ is uniformly $M$-flat if and only if $F\otimes _{S}L\leq
_{S}F\otimes _{S}M$ for every finitely generated $S$-subsemimodule $L\leq
_{S}M.$
\end{proposition}

\begin{proposition}
\label{middle-f}Let $0\longrightarrow M_{1}\overset{\gamma }{\longrightarrow
}M\overset{\delta }{\longrightarrow }M_{2}\longrightarrow 0$ be an exact
sequence of left $S$-semimodules and assume that $F_{S}$ is uniformly $M$%
-flat.

\begin{enumerate}
\item $F_{S}$ is uniformly $M_{1}$-flat.

\item If $F_{S}\in \mathcal{I}_{S}(M_{2}),$ then $F$ is uniformly $M_{2}$%
-flat.
\end{enumerate}
\end{proposition}

\begin{Beweis}
Assume that $F_{S}$ is uniformly $M$-flat.

\begin{enumerate}
\item Let $U\leq _{S}^{u}M_{1}.$ Since $M_{1}\leq _{S}^{u}M,$ we have $U\leq
_{S}^{u}M,$ whence $F\otimes _{S}U\leq _{S}^{u}F\otimes _{S}M$ and so $%
F\otimes _{S}U\leq _{S}^{u}F\otimes _{S}M_{1}$ (e.g. \cite[Lemma 3.8 (1-b)]%
{Abu}). Consequently, $F_{S}$ is uniformly $M_{1}$-flat.

\item Let $U\leq _{S}^{u}M_{2}$ and consider $\widetilde{U}:=\{m\in M\mid
\delta (m)\in U\}.$ Then $\widetilde{U}\leq _{S}^{u}M$ and we have a
commutative diagram of left $S$-semimodules with exact rows and columns%
\begin{equation*}
\xymatrix{ & & 0 \ar[d] & 0 \ar[d] & & \\ 0 \ar[r] & M_1
\ar[r]^{\tilde{\gamma}} \ar@{=}[d] & \tilde{U} \ar[r]^{\tilde{\delta}}
\ar[d]^{\tilde{\iota}} & U \ar[d]^{\iota} \ar[r] & 0 \\ 0 \ar[r] & M_1
\ar[r]^{\gamma} & M \ar[r]^{\delta} & M_2 \ar[r] & 0}
\end{equation*}%
Tensoring with $F_{S},$ we obtain a commutative diagram of Abelian monoids%
\begin{equation*}
\xymatrix{ & & 0 \ar[d] & 0 \ar@{.>}[d] & & \\ & F \otimes_S M_1
\ar[r]^{{\rm id}_F \otimes_S \tilde{\gamma}} \ar@{=}[d] & F \otimes_S
\tilde{U} \ar[r]^{{\rm id}_F \otimes_S \tilde{\delta}} \ar[d]^{{\rm id}_F
\otimes_S \tilde{\iota}} & F \otimes_S U \ar[d]^{{\rm id}_F \otimes_S \iota}
\ar[r] & 0 \\ 0 \ar[r] & F \otimes_S M_1 \ar[r]^{{\rm id}_F \otimes_S
\gamma} & F \otimes_S M \ar[r]^{{\rm id}_F \otimes_S \delta} & F \otimes_S
M_2 \ar[r] & 0}
\end{equation*}%
Since $F_{S}$ is uniformly flat, the second row is exact. By Proposition \ref%
{r-exact}, the first row is semi-exact and $\mathrm{id}_{F}\otimes _{S}%
\widetilde{\delta }$ is uniform. It follows by Lemma \ref{diagram}
\textquotedblleft 1-a\textquotedblright\ that $\mathrm{id}_{F}\otimes
_{S}\iota $ is injective. Since $F\in \mathcal{I}_{S}(M),$ we have $F\otimes
_{S}U\leq _{S}^{u}F\otimes _{S}M_{2}.$ Consequently, $F_{S}$ is uniformly $%
M_{2}$-flat.$\blacksquare $
\end{enumerate}
\end{Beweis}

Let%
\begin{equation}
A\overset{f}{\longrightarrow }B\overset{g}{\longrightarrow }C  \label{ABC-2}
\end{equation}%
be a sequence of left $S$-semimodules. The proof of the following result is
straightforward:

\begin{proposition}
\label{FP-exact}

\begin{enumerate}
\item If $F_{S}$ is a free right $S$-semimodule and \emph{(}\ref{ABC-2}\emph{%
)} is exact \emph{(}resp. semi-exact, quasi-exact, proper-exact\emph{)},
then the sequence%
\begin{equation}
F\otimes _{S}A\overset{\mathrm{id}_{F}\otimes _{S}f}{\longrightarrow }%
F\otimes _{S}B\overset{\mathrm{id}_{F}\otimes _{S}g}{\longrightarrow }%
F\otimes _{S}C  \label{Fx}
\end{equation}%
of Abelian monoids is exact \emph{(}resp. semi-exact, proper-exact,
quasi-exact\emph{)}.

\item Every free $S$-semimodule is uniformly flat.
\end{enumerate}
\end{proposition}

\begin{corollary}
\begin{enumerate}
\item If $P_{S}$ is a projective right $S$-semimodule and \emph{(}\ref{ABC-2}%
\emph{)} is exact \emph{(}resp. semi-exact, quasi-exact, proper-exact\emph{)}%
, then the sequence%
\begin{equation}
P\otimes _{S}A\overset{\mathrm{id}\otimes _{S}f}{\longrightarrow }P\otimes
_{S}B\overset{\mathrm{id}\otimes _{S}g}{\longrightarrow }P\otimes _{S}C
\label{Px}
\end{equation}%
of Abelian monoids is exact \emph{(}resp. semi-exact, proper-exact,
quasi-exact\emph{)}.

\item Every projective $S$-semimodule is uniformly flat.
\end{enumerate}
\end{corollary}

\begin{Beweis}
\begin{enumerate}
\item This follows directly from Proposition \ref{FP-exact} and Lemma \ref%
{comm} \textquotedblleft 2\textquotedblright .

\item This follows directly from the definition and \textquotedblleft
1\textquotedblright .$\blacksquare $
\end{enumerate}
\end{Beweis}

\begin{definition}
Let $Q$ be a right $S$-semimodule. We say that $Q_{S}$ (\emph{uniformly})
\emph{cogenerates} a class $\mathcal{M}$ of right semimodules iff the
following holds: for every morphisms $\iota :U\longrightarrow M$ with $M\in
\mathcal{M},$ if $(\iota ,Q):\mathrm{Hom}_{S}(M,Q)\longrightarrow \mathrm{Hom%
}_{S}(U,Q)$ is surjective (and uniform), then $0\longrightarrow U\overset{%
\iota }{\longrightarrow }M$ is injective (and uniform). If $Q$ (uniformly)
cogenerates all left $S$-semimodules, then we say that $Q$ is a (\emph{%
uniform}) \emph{cogenerator} in $\mathbb{S}_{S}.$
\end{definition}

\begin{ex}
The assumption that $S$ has an injective cogenerator might be empty. For
example the semiring $\mathbb{N}_{0}$ has no injective cogenerators.
\end{ex}

\begin{proposition}
\label{flat-inj}Let $_{T}F_{S}$ be a $(T,S)$-bisemimodule, $M$ a left $S$%
-semimodule and $X$ a left $T$-semimodule.

\begin{enumerate}
\item Let $_{T}X$ be uniformly $F\otimes _{S}M$-injective. If $F_{S}$ is
uniformly $M$-flat, then $\mathrm{Hom}_{T}(F,X)$ is uniformly $M$-injective.

\item Let $M$ be uniformly $X$-cogenerated. If $\mathrm{Hom}_{T}(F,X)$ is
uniformly $M$-injective, then $F_{S}$ is uniformly $M$-flat.
\end{enumerate}
\end{proposition}

\begin{Beweis}
Let $M$ be a left $S$-semimodule, $U\leq _{S}^{u}M$ and consider the
following commutative diagram%
\begin{equation*}
\xymatrix{\mathrm{Hom}_{S}(M,\mathrm{Hom}_{T}(F,X)) \ar[rr]^{\simeq}
\ar[dd]_{(-,\mathrm{Hom}_{T}(F,X))}& & \mathrm{Hom}_{T}(F\otimes _{S}M,X)
\ar[dd]^{(\mathrm{id}_{F}\otimes _{S}\iota _{U},X)} \\ & & \\
\mathrm{Hom}_{S}(U,\mathrm{Hom}_{T}(F,X)) \ar[rr]_{\simeq} & &
\mathrm{Hom}_{T}(F\otimes _{S}U,X)}
\end{equation*}

\begin{enumerate}
\item Let $_{T}X$ be uniformly $F\otimes _{S}M$-injective. If $F_{S}$ is
uniformly $M$-flat, then $F\otimes _{S}U\leq _{T}^{u}F\otimes _{S}M,$ whence
$(\mathrm{id}_{F}\otimes _{S}\iota _{U},X)$ is surjective and uniform.
Consequently, $(-,\mathrm{Hom}_{T}(F,X))$ is surjective and uniform. This
means that $\mathrm{Hom}_{T}(F,X)$ is uniformly $M$-injective.

\item Let $M$ be uniformly $X$-cogenerated. If $\mathrm{Hom}_{T}(F,X)$ is
uniformly $M$-injective, then $(-,\mathrm{Hom}_{T}(F,X))$ is surjective and
uniform, whence $(\mathrm{id}_{F}\otimes _{S}\iota _{U},X)$ is surjective
and uniform. So, $\mathrm{id}_{F}\otimes _{S}\iota _{U}$ is injective and
uniform. This means that $F_{S}$ is uniformly $M$-flat.$\blacksquare $
\end{enumerate}
\end{Beweis}

\begin{theorem}
\label{ic-Q}Let $_{T}F_{S}$ be a $(T,S)$-bisemimodule and assume that $_{T}%
\mathbb{S}$ has a uniformly injective-cogenerator $\mathbf{Q}.$ Then $F_{S}$
is uniformly flat if and only if $\mathrm{Hom}_{T}(F,\mathbf{Q})$ is
uniformly injective.
\end{theorem}

The analogous of Baer's criterion for injective modules over rings
\textquotedblleft $M$\emph{\ is }$R$\emph{-injective }$\Longrightarrow $%
\emph{\ }$M$\emph{\ is injective}\textquotedblright\ might fail for
semimodules over semirings.

\begin{ex}
(\cite{Ili2008}) The semifield $(\mathbb{Q}^{+},+,\cdot )$ has only two
ideals $\{0\}$ and $\mathbb{Q}^{+}$ whence every semimodule is $\mathbb{Q}%
^{+}$-injective. However, $\{0\}$ is the only injective $\mathbb{Q}^{+}$%
-semimodule (e.g. by \cite{Ili2008}).
\end{ex}

The above example motivates the following definitions:

\begin{definition}
We say that the semiring $S$ is a \emph{left} (\emph{uniformly}) \emph{%
Baer's semiring} iff every (uniformly) injective left $S$-semimodule is
(uniformly) injective. The right (uniformly) Baer-injective semirings can be
defined analogously.
\end{definition}

\begin{proposition}
Let $_{T}F_{S}$ be a $(T,S)$-bisemimodule and assume that $_{T}\mathbb{S}$
has a uniformly injective cogenerator $\mathbf{Q}.$ If $S$ is a left
uniformly Baer semiring, then the following are equivalent:

\begin{enumerate}
\item $F_{S}$ is uniformly flat;

\item For every uniform left ideal $_{S}I\leq ^{u}S,$ we have $F\otimes
_{S}I\leq _{S}^{u}F\otimes _{S}S.$
\end{enumerate}
\end{proposition}

\begin{Beweis}
We need only to prove \textquotedblleft 2\textquotedblright\ $%
\Longrightarrow $ \textquotedblleft 1\textquotedblright . Let $_{S}I\leq
_{S}^{u}S$ be a left uniform ideal. By assumption, $0\longrightarrow
F\otimes _{S}I\overset{\mathrm{id}_{F}\otimes _{S}\iota }{\longrightarrow }%
F\otimes _{S}S$ is exact and $\mathrm{id}_{F}\otimes _{S}\iota $ is a
uniform morphism of left $T$-semimodule, whence $\mathrm{Hom}_{T}(F\otimes
_{S}I,\mathbf{Q})\overset{(\mathrm{id}_{F}\otimes _{S}\iota ,\mathbf{Q})}{%
\longrightarrow }\mathrm{Hom}_{S}(F\otimes _{S}S,\mathbf{Q})\longrightarrow
0 $ is exact and $(\mathrm{id}_{F}\otimes _{S}\iota ,\mathbf{Q})$ is
uniform. Notice that $\mathrm{Hom}_{T}(F\otimes _{S}I,\mathbf{Q})\simeq
\mathrm{Hom}_{S}(I,\mathrm{Hom}_{T}(F,\mathbf{Q}))$ and $\mathrm{Hom}%
_{S}(F\otimes _{S}S,\mathbf{Q})\simeq \mathrm{Hom}_{S}(S,\mathrm{Hom}_{T}(F,%
\mathbf{Q})),$ whence $\mathrm{Hom}_{S}(I,\mathrm{Hom}_{T}(F,\mathbf{Q}))%
\overset{(\iota ,\mathrm{Hom}_{T}(F,\mathbf{Q}))}{\longrightarrow }\mathrm{%
Hom}_{S}(S,\mathrm{Hom}_{T}(F,\mathbf{Q}))\longrightarrow 0$ is exact and $%
(\iota ,\mathrm{Hom}_{T}(F,\mathbf{Q}))$ is uniform, \emph{i.e. }$\mathrm{Hom%
}_{T}(F,\mathbf{Q})$ is uniformly $S$-injective. Since $S$ is a left
uniformly Baer semiring, we conclude that $\mathrm{Hom}_{T}(F,\mathbf{Q})$
is uniformly injective as a left $S$-semimodule, whence $F_{S}$ is uniformly
flat by Theorem \ref{ic-Q}.$\blacksquare $
\end{Beweis}

\begin{theorem}
\label{lim-flat}Let $(F_{j},\{f_{jj^{\prime }}\})_{J}$ be a directed system
of right $S$-semimodules.

\begin{enumerate}
\item If each $F_{j}$ is uniformly $M$-flat, for some left $S$-semimodule $%
M, $ then $\lim\limits_{\longrightarrow }F_{j}$ is uniformly $M$-flat.

\item If each $F_{j}$ is uniformly flat, then $\lim\limits_{\longrightarrow
}F_{j}$ is uniformly flat.
\end{enumerate}
\end{theorem}

\begin{Beweis}
We need only to prove \textquotedblleft 1\textquotedblright . Assume that $%
F_{j}$ is uniformly $M$-flat for every $j\in J.$ Let $U\leq _{S}^{u}M.$ Then
$F_{j}\otimes _{S}U\leq _{S}^{u}F_{j}\otimes _{S}M$ for each $j\in J.$ It
follows by Corollary \ref{ad-l-cor} that $\lim\limits_{\longrightarrow
}(F_{j}\otimes _{S}U)\leq _{S}^{u}\lim\limits_{\longrightarrow
}(F_{j}\otimes _{S}M)$ and so we are done (note that $\lim\limits_{%
\longrightarrow }F_{j}\otimes _{S}U\simeq \lim\limits_{\longrightarrow
}(F_{j}\otimes _{S}U)$ and $\lim\limits_{\longrightarrow }F_{j}\otimes
_{S}M\simeq \lim\limits_{\longrightarrow }(F_{j}\otimes _{S}M)$).$%
\blacksquare $
\end{Beweis}

\begin{corollary}
\label{fg-sub}If every finitely generated subsemimodule of an $S$-semimodule
$F$ is uniformly flat, then $F$ is uniformly flat.
\end{corollary}

\begin{Beweis}
This follows directly from Theorem \ref{lim-flat} and the fact that every
semimodule is the direct limit of its finitely generated subsemimodules (cf.
Proposition \ref{lim-sub}).$\blacksquare $
\end{Beweis}

As a direct consequence of Theorem \ref{lim-flat} we obtain:

\begin{corollary}
Every flat $S$-semimodule is uniformly flat.
\end{corollary}

We finish this manuscript with the following open question:

\vspace*{5mm}

\textbf{Question: }When is every \emph{uniformly flat} $S$-semimodule \emph{%
flat}?

\vspace*{5mm}

\textbf{Acknowledgments.} The author thanks H. Al-Thani, Y. Katsov and A.
Patchkoria for providing him with several related papers and preprints.

\end{document}